\documentclass[12pt,a4paper,reqno]{amsart}

\usepackage{amsmath}
\usepackage{geometry}
\usepackage{amssymb, latexsym}
\usepackage{verbatim}
\usepackage{enumerate}
\usepackage{mathtools}\usepackage[tableposition=top]{caption}\usepackage{booktabs,dcolumn}
\newtheorem{theorem}{Theorem}[section]
\newtheorem{proposition}[theorem]{Proposition}

\newtheorem{lemma}[theorem]{Lemma}

\newtheorem{definition}[theorem]{Definition}

\newcommand\R{\mathbb{R}}
\newcommand\Z{\mathbb{Z}}
\newcommand\N{\mathbb{N}}
\newcommand\C{\mathbb{C}}

\newcommand\T{\mathbb{T}}
\newcommand\eps{\varepsilon}
  \newcommand\nopf{(1+\varepsilon')}

\renewcommand{\phi}{\varphi}

\begin{document}
\title[Correlations of von Mangoldt and divisor functions II]{Correlations of the von Mangoldt and higher divisor functions II.  Divisor correlations in short ranges}

\author{Kaisa Matom\"aki}
\address{Department of Mathematics and Statistics \\
University of Turku, 20014 Turku\\
Finland}
\email{ksmato@utu.fi}

\author{Maksym Radziwi{\l}{\l}}
\address{
  Department of Mathematics\\
  Caltech\\
  1200 E California Blvd\\
  Pasadena, CA, 91125\\
}
\email{maksym.radziwill@gmail.com}

\author{Terence Tao}
\address{Department of Mathematics, UCLA\\
405 Hilgard Ave\\
Los Angeles CA 90095\\
USA}
\email{tao@math.ucla.edu}
\begin{abstract}  We study the problem of obtaining asymptotic formulas for the sums $\sum_{X < n \leq 2X} d_k(n) d_l(n+h)$ and $\sum_{X < n \leq 2X} \Lambda(n) d_k(n+h)$, where $\Lambda$ is the von Mangoldt function, $d_k$ is the $k^{\operatorname{th}}$ divisor function, $X$ is large and $k \geq l \geq 2$ are integers. We show that for almost all $h \in [-H, H]$ with $H = (\log X)^{10000 k \log k}$, the expected asymptotic estimate holds. In our previous paper we were able to deal also with the case of $\Lambda(n) \Lambda(n + h)$ and we obtained better estimates for the error terms at the price of having to take $H = X^{8/33 + \varepsilon}$. 

\end{abstract}

\maketitle

\section{Introduction}

This paper is a sequel to our previous work \cite{mrt-corr}.  As in our previous paper, we are interested in the correlations\footnote{The results here can also be applied to the sums $\sum_{n \leq X} d_k(n) d_l(n+h)$ and $\sum_{n \leq X} \Lambda(n) d_k(n+h)$ after some minor modifications to the coefficients of the polynomials $P_{k,l,h}$ and $Q_{k,h}$; we leave the details to the interested reader.}
\begin{align} \label{correlationformulas}
\sum_{X < n \leq 2X} d_k(n) d_l(n+h) \quad \text{and} \quad
\sum_{X < n \leq 2X} \Lambda(n) d_k(n+h),
\end{align}
where $\Lambda$ is the von Mangoldt function and  for $k$ integer, 
$d_k(n) \coloneqq  \sum_{n_1 \dotsm n_k = n} 1$
is the $k^{\operatorname{th}}$ divisor function. 
In contrast to \cite{mrt-corr}, we omit the case of $\Lambda(n)\Lambda(n + h)$ from consideration as we will have nothing new to say for these correlations. 

When $h \neq 0$ is fixed, $X$ goes to infinity, and $k \geq l \geq 2$ are fixed integers, there are well-established conjectures for the asymptotic values of expressions in \eqref{correlationformulas}. For instance it is conjectured in \cite{vinogradov}, \cite{ivic}, \cite[Conjecture 3]{conrey} that
\begin{equation} \label{dkdk}
\sum_{X < n \leq 2X} d_k(n) d_l(n + h) = P_{k,l,h}(\log X) \cdot X + O_{\varepsilon}(X^{1/2 + \varepsilon})
\end{equation}
and
\begin{equation} \label{ldk}
\sum_{X < n \leq 2X} \Lambda(n) d_k(n + h) = Q_{k,h}(\log X) X + O_{\varepsilon}(X^{1/2 + \varepsilon})
\end{equation}
with $P_{k,l,h}$, $Q_{k,h}$ explicitly computable polynomials of degree $k+l-2$ and $k - 1$ respectively for $h \neq 0$. This is also wide open as soon as $k,l \geq 3$. 
See \cite{brudern},\cite{mrt-corr} for further discussion of these conjecture as well as previous progress on them.

In \cite{mrt-corr} the conjectures \eqref{dkdk} and \eqref{ldk} were proven on average for almost all shifts $h$ in a short interval, with weaker error terms and for integer $k,l > 0$. Precisely, for $H = X^{8/33 + \varepsilon}$ we showed that $\eqref{dkdk}$ and \eqref{ldk} hold with an error term that is $\ll X (\log X)^{-A}$ for all but at most $\ll_{A} X (\log X)^{-A}$ shifts $h \in [1, H]$. This improved on earlier results which contained the case $H = X^{1/3 + \varepsilon}$ either in explicit or implicit form (see \cite{mikawa}, \cite{bbmz}). 

In this paper we are interested in substantially smaller $H$'s, specifically $H$ as small as $\log^B X$ for some large constant $B>0$. 
Our error terms are considerably weaker, therefore we are only able to confirm on average the leading term of the polynomials $P_{k,l,h}, Q_{k,h}$ appearing in conjectures \eqref{dkdk} and \eqref{ldk}.

\begin{theorem}[Main theorem]\label{uncond} Let $0 < \eps < 1/2$ and $k \geq l \geq 2$ be fixed integers.  Suppose that $\log^{10000 k \log k} X \leq H \leq X^{1-\eps}$ for some $X \geq 2$. There exists constants $C_{k,l,h} > 0$ and $C_{k,h} > 0$ such that as $X \rightarrow \infty$, \begin{itemize} \item [(i)] 
  \begin{align*}
\sum_{0 < |h| \leq H} \Big | \sum_{X < n \leq 2X} & d_k(n) d_{l}(n + h) - C_{k,l,h} \cdot X (\log X)^{k + l - 2} \Big | = o \Big ( H X (\log X)^{k + l - 2} \Big )
  \end{align*}
\item[(ii)]
  $$ \sum_{0 < |h| \leq H} \Big | \sum_{X < n \leq 2X} \Lambda(n) d_k(n+h) - C_{k,h} \cdot X (\log X)^{k - 1} \Big | = o(H X (\log X)^{k - 1}).$$
\end{itemize}
In fact, the constants $C_{k,l,h}$ and $C_{k,h}$ are explicitly given by the formulae 
  $$
 C_{k,l,h} \coloneqq \frac{1}{\Gamma(k) \Gamma(l)} \prod_{p} \mathfrak{S}_{k,l,p}(h) \ , \   C_{k,h} \coloneqq  \frac{1}{\Gamma(k)} \Big ( \prod_{p} \mathfrak{S}_{k,p}(h) \Big )
  $$
  with the singular series $\mathfrak{S}_{k,l,p}(h)$ and $\mathfrak{S}_{k,p}(h)$ defined in \cite{mrt-corr}.
\end{theorem}

The saving implicit in the $o(\cdot)$ is quite small; it is of the form $(\log X)^{-\varepsilon(X)}$ for a function $\varepsilon(X)$ going to zero arbitrarily slowly.

The weakness of the error term is largely due to us mostly exploiting the anatomy of integers, in particular the multiplicativity of $d_k(n)$, and not much else. Obtaining an error term of the form $O_{A}(X \log^{-A} X)$ is an important open problem, as it likely would allow us to obtain the Hardy-Littlewood conjectures on average with a small averaging of the form $0 < |h| \leq X^{\varepsilon}$. The weakness of the error term also means that our results have no new consequences for moments of the Riemann zeta-function, where at least a power-saving would be required. 

Our methods are not able to reduce the size of $H$ below $(\log X)^{k \log k - k + 1}$, since they depend crucially on having
\begin{equation} \label{asympt}
\sum_{|h| \leq H} d_k(n + h) \sim c_{k} H (\log X)^{k - 1} . 
\end{equation}
for most integers $X < n < 2X$, and with $c_k > 0$ a constant depending only on $k$. When $H < (\log X)^{k \log k - k + 1 - \varepsilon}$ the asymptotic \eqref{asympt} is false for almost all integers $X < n < 2X$, because the main contribution to the average size of $d_k(n)$ comes from integers with $(1 + o(1)) k \log \log X$ prime factors, and such integers do not typically occur in intervals of length less than $(\log X)^{k \log k - k + 1 -\varepsilon}$. 

With some effort it is possible to generalize our results to a wide class of multiplicative functions. The most immediate extension is to $k, l \geq 2$ non-integer real numbers. This requires only an extension of Lemma \ref{majorarcmain} to $k, l \geq 2$ non-integer and follows (after some work) from the Selberg-Delange method. A bit more work should allow also to extend the result to all $k, \ell > 0$ real. In principle one should be able to formulate a rather general result that works for arbitrary multiplicative functions satisfying some minor necessary conditions, such as for example not being pretentious to characters $\chi(n) n^{it}$ of very low conductor (smaller than $(\log X)^{A}$). Doing so in this paper would however obscur an already complicated proof. 

\subsection{Overview of proof}

For sake of exposition let us restrict attention to the case $k=l$.
The first step in our proof consists (as in \cite{mr}) in throwing away from 
$$
\sum_{X < n \leq 2X} d_k(n) d_k(n + h)
$$
a small set of integers for which either $n$ or $n + h$ has an unusual factorization (e.g primes). Specifically we require that the remaining integers have
\begin{enumerate}
\item [(A)] prime factors in certain wide intervals; and
\item [(B)] not too many prime factors.
\end{enumerate}
It is already at this step that we forfeit the possibility of a respectable saving in the error term. 
  Let us use  $f_k$ to denote the restriction of the $k^{\operatorname{th}}$ divisor function $d_k$ to a set ${\mathcal S}_{k,X}$ of integers satisfying properties (A) and (B); see Section \ref{circle-sec} for a precise description of this set. Naively proceeding as in our previous paper \cite{mrt-corr}, the next step would then consist in bounding non-trivially (i.e., beyond the bounds obtained from Parseval) the quantity
  \begin{equation} \label{minorarc}
\sup_{\alpha} \int_{\alpha - 1/H}^{\alpha + 1/H} |S_{f_k}(\beta;X)|^2 d \beta , \quad \text{where} \quad S_{f_k}(\alpha;X) \coloneqq \sum_{X < n \leq 2X} f_k(n) e(n \alpha) 
\end{equation}
and $\alpha$ is restricted to be ``minor arc'' in a suitable sense. However, in contrast to our previous paper, as $H$ is so small we must use techniques coming from the work of the first two authors \cite{mr}, and so we no longer necessarily have large savings in our bounds for \eqref{minorarc}. In fact, since $f_k$ fluctuates quite wildly in size (in particular its $L^2$ norm is substantially larger than its $L^1$ norm) to succeed with the approach from \cite{mrt-corr}, one would need to save in \eqref{minorarc} an arbitrary power of the logarithm compared to the trivial bound coming from Parseval. This is possible if $H > (\log X)^{\psi(X)}$ for some $\psi(X)$ going to infinity arbitrarily slowly, after exploiting  property (A) of our set of integers. The main input then is the work of the first two authors \cite{mr} adapted to the setting of the divisor function. The main innovation compared to \cite{mr} is the use of a mean-value theorem that plays well with the wild fluctuation of the size of $d_k(n)$.

Going below to the level of $H = (\log X)^{B}$ for a large constant $B$ requires additional ideas, in particular a much subtler use of harmonic analysis. Since in the range $H = (\log X)^{B}$ we cannot obtain a saving of an arbitrary power of the logarithm in \eqref{minorarc}, we modify the argument so that it can accept any non-trivial saving in \eqref{minorarc}, at the cost of an additional input, namely a non-trivial bound for the quantity 
\begin{equation} \label{largevalues}
\sum_{j = 1}^{J} \int_{\alpha_j - 1/H}^{\alpha_j + 1/H} |S_{f_k}(\beta;X)|^2 d \beta 
\end{equation}
whenever the $\alpha_j$ are $1/H$-separated, where by ``non-trivial'' we mean that the bound for \eqref{largevalues} improves over the trivial bound obtained from bounding each term individually as soon as $J \rightarrow \infty$.
It is in the proof of such a large value estimate that the property (B) of our set of integers is used crucially and repeatedly.  This property allows us to construct an efficient sieve majorant for $d_k$ on the set of integers satisfying (B).

The estimate \eqref{largevalues} is easy when $J$ is large and when $J$ is very small, so we can at the outset assume that $J$ is a small power of $H$, which corresponds to the difficult middle-range. As a first step towards establishing \eqref{largevalues} we use Gallagher's lemma to reduce the problem to the estimation of
$$
\sum_{j = 1}^{J} \int_{X}^{2X} \left| \sum_{x \leq n \leq x + H} f_k(n) e(\alpha_j n) \right|^2\ d x 
$$
Dualizing (and doing some additional technical preparation), the problem boils down to controlling for most $x \in [X, 2X]$, 
\begin{equation} \label{dualized}
\sum_{x \leq n \leq x + H} \widetilde{d}_k(n) e(n (\alpha_j - \alpha_{i}))
\end{equation}
on average over $\alpha_j$'s that are $1/H$ separated and where $\widetilde{d}_k(n)$ is a sieve majorant for $f_k(n)$.
In \eqref{dualized} we average over $J^2$ tuples $(\alpha_j, \alpha_i)$ and we must win over the individual bounds as soon as $J \rightarrow \infty$. 
For a parameter $Q$ to be chosen later, let us use $q_{\alpha, Q}$ to denote the smallest integer $q$ for which there exists an $(a,q) = 1$ with 
$$
\Big \| \alpha - \frac{a}{q} \Big \|_{\mathbb{T}} \leq \frac{1}{q Q}. 
$$
where $\| x \|_{\mathbb{T}}$ denotes the distance of $x$ to the integers. 
If we choose $Q$ to be a small power of $H$, but a large power of $J$, then the technicalities that go before dualizing allow us to throw away from consideration a pathological set of $x$'s and to assume that in \eqref{dualized} we are saving $1/q_{\alpha_j - \alpha_i, Q}$ over the ``trivial bound'' of $H (\log X)^{k - 1}$ (in fact this ``trivial bound'' holding on the set of non-exceptional $x$'s is also the result of a considerable amount of work!). Now if the $\alpha_j$ were all $1/Q$ separated, it would be enough to split into two cases according to the size of $q_{\alpha_j - \alpha_i, Q}$. We could dispose of those $(\alpha_j, \alpha_i)$ for which $q_{\alpha_j - \alpha_i, Q} > J^2$ by simply bounding \eqref{dualized} by $(1/J)^2 H (\log X)^{k - 1}$, while for the remaining tuples $(\alpha_j, \alpha_i)$ we could show that the contribution of the other $(\alpha_j, \alpha_i)$ is negligible by using an estimate
\begin{equation} \label{greentaoestimate}
\sum_{(\alpha_j, \alpha_i) : q_{\alpha_j - \alpha_i, Q} \leq J^2} \frac{1}{q_{\alpha_j - \alpha_i, Q}} \ll_{\varepsilon} J^{1 + \varepsilon}
\end{equation}
of Green-Tao \cite{gt-selberg}, which is available under this separation hypothesis on the $\alpha_j$.  This would be sufficient, since in this context the trivial bound corresponds to $J^2$, so a bound of $J^{1 + \varepsilon}$ means that there are only few tuples for which $q_{\alpha_j - \alpha_i, Q}$ is very small.  
The main input in \eqref{greentaoestimate} is the large sieve.

However the situation is a bit more complicated because the $\alpha_j$'s are $1/H$ spaced and not $1/Q$ spaced. Since $1/H$ is significantly smaller than $1/Q$ this means that we could have $q_{\alpha_j - \alpha_i, Q} = 1$ for all tuples $(\alpha_j, \alpha_i)$. In such a situation the left-hand side of \eqref{greentaoestimate} would be $\gg J^2$, and now this estimate fails.
To circumvent this issue, we group the frequencies $\alpha_j$ into ``batches'' of $1/Q$-nearby frequencies (i.e., the set of the frequencies within the same batch has diameter at most $1/Q$). We pick a representative frequency $\beta_i$ from each batch $\mathcal{I}_i$ so that the $\beta_i$ are roughly $1/Q$ separated, and there are $|\mathcal{I}_{i}|$ frequencies $\alpha_j$ in the batch associated to $\beta_i$. We then perform the same argument as before, but for the batches (for example the analogue of the condition $q_{\beta_i - \beta_j, Q} > J^2$ for batches is that we have $q_{\alpha_{i'} - \alpha_{j'}, Q} > J^2$ for all $\alpha_{i'} \in \mathcal{I}_{i}$ and all $\alpha_{j'} \in \mathcal{I}_j$). The main technical difference is that we need a generalization of the result of Green-Tao to sums of the form,
\begin{equation} \label{generalizedsum}
\sum_{(\beta_j, \beta_i) : q_{\beta_j - \beta_i,Q} \leq J^2} \frac{|\mathcal{I}_i|^{1/2} \cdot |\mathcal{I}_{j}|^{1/2}}{q_{\beta_i - \beta_j, Q}} .
\end{equation}
There are however no difficulties involved with this generalization, and we bound \eqref{generalizedsum} by $J^{5/4}$ which is sufficient (and with a bit more work we could have improved this bound to $O_{\varepsilon}(J^{1 + \varepsilon})$).

\subsection{Acknowledgments}

KM was supported by Academy of Finland grant no. 285894. MR was supported by a NSERC DG grant, the CRC program and a Sloan Fellowship. TT was supported by a Simons Investigator grant, the James and Carol Collins Chair, the Mathematical Analysis \& Application Research Fund Endowment, and by NSF grant DMS-1266164.

Part of this paper was written while the authors were in residence at MSRI in Spring 2017, which is supported by NSF grant DMS-1440140.

The authors are grateful to Jiseong Kim for pointing out a few typos that appeared in the published version of this paper.

\section{Notation and preliminaries}\label{notation-sec}

All sums and products will be over integers unless otherwise specified, with the exception of sums and products over the variable $p$ (or $p_1$, $p_2$, $p'$, etc.) which will be over primes.  To accommodate this convention, we adopt the further convention that all functions on the natural numbers are automatically extended by zero to the rest of the integers, e.g. $\Lambda(n) = 0$ for $n \leq 0$.

We use $A = O(B)$, $A \ll B$, or $B \gg A$ to denote the bound $|A| \leq C B$ for some constant $C$.  If we permit $C$ to depend on additional parameters then we will indicate this by subscripts, thus for instance $A = O_{k,\eps}(B)$ or $A \ll_{k,\eps} B$ denotes the bound $|A| \leq C_{k,\eps} B$ for some $C_{k,\eps}$ depending on $k,\eps$.  If $A,B$ both depend on some large parameter $X$, we say that $A = o(B)$ as $X \to \infty$ if one has $|A| \leq c(X) B$ for some function $c(X)$ of $X$ (as well as further ``fixed'' parameters not depending on $X$), which goes to zero as $X \to \infty$ (holding all ``fixed'' parameters constant). We also write $A \asymp B$ for $A \ll B \ll A$, with the same subscripting conventions as before.

We use $ \T \coloneqq  \R/\Z$ to denote the unit circle, and $e: \T \to \C$ to denote the fundamental character
$$ e(x) \coloneqq  e^{2\pi i x}.$$

We use $1_E$ to denote the indicator of a set $E$, thus $1_E(n) = 1$ when $n \in E$ and $1_E(n) = 0$ otherwise. Similarly, if $S$ is a statement, we let $1_S$ denote the number $1$ when $S$ is true and $0$ when $S$ is false, thus for instance $1_E(n) = 1_{n \in E}$.  If $E$ is a finite set, we use $\# E$ to denote its cardinality.

We use $(a,b)$ and $[a,b]$ for the greatest common divisor and least common multiple of natural numbers $a,b$ respectively, and write $a|b$ if $a$ divides $b$.  We also write $a = b\ (q)$ if $a$ and $b$ have the same residue modulo $q$.

Given a sequence $f\colon \N \to \C$, we define the $\ell^2$ norm $\|f\|_{\ell^2}$ of $f$ as
$$ \|f\|_{\ell^2} \coloneqq  \left(\sum_n |f(n)|^2\right)^{1/2}$$
and similarly define the $\ell^\infty$ norm 
$$ \|f\|_{\ell^\infty} \coloneqq  \sup_n |f(n)|.$$
If $f$ is finitely supported, we define the exponential sum $S_f \colon \R/\Z \to \C$ by the formula
$$ S_f(\alpha;X) \coloneqq \sum_{X < n \leq 2X} f(n) e(n\alpha).$$

Given two arithmetic functions $f, g \colon \N \to \C$, the Dirichlet convolution $f \ast g$ is defined by
$$ f \ast g(n) \coloneqq  \sum_{d|n} f(d) g\left(\frac{n}{d}\right).$$

Given an interval $I$, we set
$$
\omega(n ; I ) = \sum_{\substack{p | n \\ p \in I}} 1
$$
and
$$
\Omega(n ; I) = \sum_{\substack{p^{\alpha} || n \\ p \in I}} \alpha
$$
where
$p^{\alpha} || n$ means that $p^{\alpha} | n$ but $p^{\alpha + 1} {\not|} n$.

\subsection{Fourier Analysis}

We shall frequently use the following Fourier analytic observation of Gallagher \cite[Lemma 1]{gallagher}.
\begin{lemma}
\label{le:Gallagher}
Let $f\colon \mathbb{N} \to \mathbb{R}$ be finitely supported, and let $Y \geq 1$. Then
\[
\int_{-\frac{1}{2Y}}^{\frac{1}{2Y}} |S_f(\alpha;X)|^2 d\alpha \ll \int_{-\infty}^\infty \left|\frac{1}{Y} \sum_{\substack{x \leq n \leq x+Y \\ X < n \leq 2X}} f(n)\right|^2 dx.
\]
\end{lemma}

\subsection{Divisor bounds}

We will need several standard bounds on the divisor functions $d_k$ and related multiplicative functions. The standard divisor bound gives
\begin{equation}\label{divisor-bound}
 \sum_{1 \leq n \leq x} d_k(n)^l \ll_{k,l} x \log^{k^l-1} x,
\end{equation}
for any $k,l \geq 1$ and $x \geq 2$, which can be obtained by elementary number theory methods (see e.g., \cite[formula (1.80)]{ik}).  In particular, we have the crude divisor bound
\begin{equation}\label{crude}
 d_k(n) \ll n^{o(1)}
\end{equation}
for fixed $k$ and $n \to \infty$.  

The exponent of $\log x$ in the right-hand side of \eqref{divisor-bound} is often too large for our desired applications when $l \geq 2$. To partially circumvent this issue, we will rely on some consequences of the upper bounds of Henriot~\cite{henriot, henriot2}, which do not lose excess powers of $\log X$ as long as there are non-trivial shifts in the factors of $d_k$ on the left-hand side. But let us first define a class of multiplicative functions to which such results apply.

\begin{definition}
For $A, B \geq 1$ and $\varepsilon > 0$, let $\mathcal{M}(A, B, \varepsilon)$ denote the class of multiplicative functions $f \colon \mathbb{N} \to \mathbb{R}_{\geq 0}$ such that
\begin{itemize}
\item[(i)] $f(p^\nu) \leq A^\nu$ \text{ for all $p \in \mathbb{P}$ and $\nu \geq 1$;}
\item[(ii)] $f(n) \leq B n^\varepsilon$ \text{ for all $n \geq 1$.}
\end{itemize}
\end{definition}
Now we are ready to state bounds for multiplicative functions and their correlations.

\begin{lemma}\label{henriot}  Let $\delta \in (0, 1)$ be fixed and let $X \geq Y \geq X^{\delta} \geq 2$.
\begin{itemize}
\item[(i)] Let $A, B \geq 1$ and let $f \in \mathcal{M}(A, B, \delta^3/1000)$. Then
\begin{equation}
\label{eq:ShiuBoundf2}
\frac{1}{Y} \sum_{X < n \leq X+Y} f(n) \ll_{\delta, A, B} \prod_{p \leq X} \left(1+\frac{f(p)-1}{p}\right).
\end{equation}
\item[(ii)] Let $A \geq 1$ and let $k, l \in \mathbb{N}$. If $1 \leq q \leq \log^{A} X$, and $h_1,\dots,h_l \ll X$ are distinct integers, then
$$ \sum_{X < n \leq X + Y} d_k(qn+h_1) \dotsm d_k(qn+h_l) \ll_{\delta, A,k,l} Y \log^{l(k-1)+o(1)} X.$$
\item[(iii)] Let $k, m, r \in \mathbb{N}$. Let $L_j(h)$ ($j = 1, \dotsc, m$) be distinct linear forms of $h = (h_1, \dotsc, h_r)$ with integer coefficients bounded in magnitude by $C$, and let $H \geq 1$. Then
\[
\sum_{\substack{h_1, \dotsc, h_r \\ |h_j| \leq H \\ L_i(h) \neq L_j(h) \text{ for all $i \neq j$}}} \sum_{X < n \leq X + Y} d_k(n + L_1(h)) \dotsm d_k(n + L_m(h)) \ll_{\delta, k, m, r, C} H^r Y (\log X)^{m(k-1)}.
\]
\item[(iv)] Let $A, B \geq 1$ and let $f_1, f_2 \in \mathcal{M}(A, B, \delta^3/1000)$. Let $X^\delta \geq q \geq 1$ and $H \geq 1$. Then
\[
\begin{split}
&
\sum_{0 < |h| \leq H} \sum_{\substack{X < n \leq X+Y \\ (n, q) = 1}} f_1(n) f_2(n+hq)\\
&\quad  \ll_{
\delta, A, B} HY\prod_{\substack{p \leq X  \\ p \nmid q}} \left(1+\frac{f_1(p)-1}{p}\right) \left(1+\frac{f_2(p)-1}{p}\right) \prod_{p \mid q} \left(1-\frac{1}{p}\right).
\end{split}
\]
\end{itemize}
\end{lemma}

\begin{proof}
(i) This is essentially due to Shiu~\cite{shiu} but in exactly this form follows from Henriot's work~\cite[Theorem 3]{henriot}. We shall give details of proofs of claims (ii) and (iii) which are slightly more complicated applications of Henriot's result, and leave the details of (i) and (iv) for the interested reader.

(ii) Let us write, for each $j = 1, \dotsc, l$, $q_j = q/(q, h_j)$ and $a_j = h_j/(q, h_j)$. Then $(a_j, q_j) = 1$ for every $j$, and additionally for $i \neq j$ we have either $q_i \neq q_j$ or $a_i \neq a_j$. Furthermore, since $q \leq (\log X)^A$, we have $d_k((q, h_j)) = (\log X)^{o(1)}$ for each $j$. Hence it suffices to show that
$$ \sum_{X < n \leq X + Y} d_k(q_1 n+a_1) \dots d_k(q_l n+a_l) \ll_{\delta, A,k,l} X \log^{l(k-1)+o(1)} X$$
where $a_j$ and $q_j$ satisfy the conditions above.

In the notation of Henriot's paper~\cite{henriot} we have 
\[
\begin{split}
Q_j(n) &= q_j n + a_j, \quad D = \prod_{i < j} (q_j a_i - a_j q_i)^2, \quad F(n_1, \dotsc, n_l) = d_k(n_1) \dotsm d_k(n_l), \\
\rho_{Q_j}(n) &= \#\{u \pmod{n} \colon q_j u+a_j \equiv 0 \pmod{n}\} =
\begin{cases}
1, & \text{if $(n, q_j) = 1$;} \\
0 & \text{otherwise,}
\end{cases} \\
\rho(p) &= \#\{n \pmod{p} \colon \prod_{j=1}^l (q_j n + a_j) \equiv 0 \pmod{p}\} = 
\begin{cases}
l & \text{if $p \nmid q_1 \dotsc q_l D$;}\\
\leq l & \text{otherwise,}
\end{cases}
\end{split}
\]
and
\[
\begin{split}
\Delta_D &= \prod_{p \mid D} \left(1 + \sum_{\emptyset \neq J \subset \{1, \dotsc, l\}} k^{|J|} \frac{\# \{n \pmod{p^2} \colon p \mid q_j n + a_j \text{ for all $j \in J$}\} }{p^2} \right) \\
&= \prod_{p \mid D} \left(1+\frac{O_{k, l}(1)}{p}\right).
\end{split}
\]
Writing $D' = D q_1 \dotsm q_l$, we get from~\cite[Theorem 3]{henriot} that
\[
\begin{split}
&\frac{1}{Y} \sum_{X < n \leq X + Y} d_k(q_1 n+a_1) \dotsm d_k(q_l n+a_l) \\
&\ll_{\delta, k, l} \Delta_D \prod_{\substack{p \leq X}} \left(1-\frac{\rho(p)}{p}\right) \sum_{\substack{n_1 \dotsm n_l \leq X \\ (n_1 \dotsm n_l, D) = 1}} d_k(n_1) \dotsm d_k(n_l) \cdot \frac{\rho_{Q_1}(n_1) \dotsc \rho_{Q_l}(n_l)}{n_1 \dotsc n_l} \\
&\ll \prod_{p \mid D} \left(1+\frac{O_{k, l}(1)}{p}\right) \prod_{\substack{p \leq X \\ p \nmid D'}} \left(1-\frac{l}{p}\right) \sum_{\substack{n_1 \dotsm n_l \leq X \\ (n_1 \dotsm n_l, D') = 1}} \frac{d_k(n_1) \dotsm d_k(n_l)}{n_1 \dotsc n_l} \\
&\ll_{A, k, l} (\log \log X)^{O_{k,l}(1)} \prod_{\substack{p \leq X \\ p \nmid D'}} \left(1-\frac{l}{p}\right) \prod_{\substack{p \leq X \\ p \nmid D'}} \left(1+\frac{k}{p}\right)^l \\
&\ll_{k,l} \log^{l(k-1)+o(1)} X.
\end{split}
\]

(iii) This time we apply Henriot's result with
\[
Q_j(n) = n + L_j(h), \quad D = \prod_{i < j} (L_i(h) - L_j(h))^2, \quad F(n_1, \dotsc, n_m) = d_k(n_1) \dotsm d_k(n_m), 
\]
\[
\rho_{Q_j}(n) = \#\{u \pmod{n} \colon u+L_j(h) \equiv 0 \pmod{n}\} = 1,
\]
\[
\rho(p) = \#\{n \pmod{p} \colon \prod_{j=1}^m (n + L_j(h)) \equiv 0 \pmod{p}\} = 
\begin{cases}
m & \text{if $p \nmid D$;}\\
\leq m & \text{otherwise,}
\end{cases}
\]
and
\[
\begin{split}
\Delta_D &= \prod_{p \mid D} \left(1 + \sum_{\emptyset \neq J \subset \{1, \dotsc, m\}} k^{|J|} \frac{\# \{n \pmod{p^2} \colon p \mid n + L_j(h) \text{ for all $j \in J$}\} }{p^2} \right) \\
&= \prod_{p \mid D} \left(1+\frac{O_{k}(1)}{p}\right).
\end{split}
\]
We get from~\cite[Theorem 3]{henriot} that
\[
\begin{split}
&\frac{1}{Y} \sum_{X < n \leq X + Y} d_k(n+L_1(h)) \dotsm d_k(n+L_m(h)) \\
&\ll_\delta \Delta_D \prod_{\substack{p \leq X}} \left(1-\frac{\rho(p)}{p}\right) \sum_{\substack{n_1 \dotsm n_m \leq X \\ (n_1 \dotsm n_m, D) = 1}} d_k(n_1) \dotsm d_k(n_m) \cdot \frac{\rho_{Q_1}(n_1) \dotsm \rho_{Q_m}(n_m)}{n_1 \dotsm n_m} \\
&\ll \prod_{p \mid D} \left(1+\frac{O_{k,m}(1)}{p}\right) \prod_{\substack{p \leq X \\ p \nmid D}} \left(1-\frac{m}{p}\right) \sum_{\substack{n_1 \dotsm n_m \leq X \\ (n_1 \dotsm n_m, D) = 1}} \frac{d_k(n_1) \dotsm d_k(n_m)}{n_1 \dotsc n_m} \\
&\ll \prod_{p \mid D} \left(1+\frac{O_{k,m}(1)}{p}\right) \prod_{p \leq X} \left(1+\frac{k-1}{p}\right)^m.
\end{split}
\]
Here by the arithmetic-geometric inequality
\[
\begin{split}
&
\prod_{p \mid D} \left(1+\frac{O_{k,m}(1)}{p}\right)\\
&\quad  \leq \sum_{i < j} \prod_{p \mid L_i(h) - L_j(h)} \left(1+\frac{O_{k,m}(1)}{p}\right)^{m(m-1)/2} \\
&\quad \leq \sum_{i < j} \prod_{p \mid L_i(h) - L_j(h)} \left(1+\frac{O_{k,m}(1)}{p}\right).
\end{split}
\]
Now
\[
\begin{split}
&\sum_{i < j} \sum_{\substack{h_1, \dotsc, h_r \\ |h_j| \leq H \\ L_i(h) \neq L_j(h) \text{for all $i \neq j$}}} \prod_{p \mid L_i(h) - L_j(h)} \left(1+\frac{O_{k,m}(1)}{p}\right) \\
&\leq H^{r-1} \sum_{\substack{h_1 = O_r(H)}} \prod_{p \mid h_1} \left(1+\frac{O_{k,m}(1)}{p}\right) \\
&\ll_r H^{r} \prod_{p = O_r(H)} \left(1+\frac{\left(1+\frac{O_{k.m}(1)}{p}\right)-1}{p}\right)  \ll_{k, m, r} H^r.
\end{split}
\]
by Shiu's bound (claim (i)).
\end{proof}

\section{Applying the circle method}\label{circle-sec}

Let $f,g: \Z \to \C$ be functions supported on a finite set, and let $h$ be an integer.  Following the Hardy-Littlewood circle method, we can express the correlation $\sum_n f(n) \overline{g}(n+h) $ as an integral
$$ \sum_n f(n) \overline{g}(n+h) = \int_\T S_f(\alpha;X) \overline{S_g(\alpha;X)} e(\alpha h)\ d\alpha.
$$
If we then designate some (measurable) portion ${\mathfrak M}$ of the unit circle $\T$ to be the ``major arcs'', we thus have
\begin{equation}\label{sam}
\sum_n f(n) \overline{g}(n+h) -  \operatorname{MT}_{{\mathfrak M},h} = \int_{\mathfrak m} S_f(\alpha;X) \overline{S_g(\alpha;X)} e(\alpha h)\ d\alpha
\end{equation}
where $ \operatorname{MT}_{{\mathfrak M},h} $ is the \emph{main term}
\begin{equation}\label{mt-def}
 \operatorname{MT}_{{\mathfrak M},h} \coloneqq  \int_{\mathfrak M} S_f(\alpha;X) \overline{S_g(\alpha;X)} e(\alpha h)\ d\alpha
\end{equation}
and ${\mathfrak m} \coloneqq  \T \backslash {\mathfrak M}$ denotes the complementary \emph{minor arcs}.

In our previous paper \cite{mrt-corr}, the following criterion was established to obtain asymptotics for such correlations on average:

\begin{lemma}\label{chebyshev}  Let $H \geq 1$ and $\eta, F, G, X > 0$.  Let $f,g: \Z \to \C$ be functions supported on a finite set, let ${\mathfrak M}$ be a  measurable subset of $\T$, and let ${\mathfrak m} \coloneqq  \T \backslash {\mathfrak M}$.  For each $h$, let $\operatorname{MT}_h$ be a complex number.  Let $h_0$ be an integer. Assume the following axioms:
\begin{itemize}
\item[(i)] (Size bounds)  One has $\|f\|_{\ell^2}^2 \ll F^2 X$ and $\|g\|_{\ell^2}^2 \ll G^2 X$.
\item[(ii)]  (Major arc estimate)  For all but $O(\eta H)$ integers $h$ with $|h-h_0| \leq H$, one has
$$ \int_{\mathfrak M} S_f(\alpha;X) \overline{S_g(\alpha;X)} e(\alpha h)\ d\alpha = \operatorname{MT}_{h} + O( \eta F G X ).$$
\item[(iii)]   (Minor arc estimate)  For each $\alpha \in {\mathfrak m}$, one has
\begin{equation}\label{fF}
 \int_{{\mathfrak m} \cap [\alpha-1/2H, \alpha+1/2H]} |S_f(\beta;X)|^2\ d\beta \ll \eta^6 F^2 X.
\end{equation}
\end{itemize}
Then for all but $O(\eta H)$ integers $h$ with $|h-h_0| \leq H$, one has
\begin{equation}\label{fg}
 \sum_n f(n) \overline{g(n+h)} = \operatorname{MT}_{h} + O( \eta FG X ).
\end{equation}
\end{lemma}

\begin{proof}  See \cite[Corollary 3.2]{mrt-corr}.
\end{proof}

This lemma is not suitable for establishing Theorem \ref{uncond}, because the $\ell^2$ norm of $d_k$ or $d_l$ is somewhat large, causing the $F$ and $G$ quantities to contain some unwanted powers of $\log X$ that will dominate the rather small gain $\eta$ that can be obtained in this case. We will therefore need the following more complicated variant of Lemma \ref{chebyshev} that avoids any direct estimation of $\ell^2$ norms of $f$ or $g$, replacing such estimates with ``large values estimates'' for $S_f$ and $S_g$.

\begin{lemma}\label{chebyshev-2}  Let $H \geq 1$, $K \geq 1$ and $\eta, F, G > 0$, and let $0 < \delta \leq 1/2$.  Let $f,g: \Z \to \C$ be functions supported on a finite set, let ${\mathfrak M}$ be a  measurable subset of $\T$, and let ${\mathfrak m} \coloneqq  \T \backslash {\mathfrak M}$.  For each $h$, let $\operatorname{MT}_h$ be a complex number.  Let $h_0$ be an integer. Assume the following axioms:
\begin{itemize}
\item[(i)]  (Major arc estimate)  For all but $O(\eta H)$ integers $h$ with $|h-h_0| \leq H$, one has
\begin{equation}\label{flip}
 \int_{\mathfrak M} S_f(\alpha;X) \overline{S_g(\alpha;X)} e(\alpha h)\ d\alpha = \operatorname{MT}_{h} + O( \eta F G X ).
\end{equation}
\item[(ii)]   (Minor arc estimate)  For each $\alpha \in {\mathfrak m}$, one has
\begin{equation}\label{fF2}
 \int_{{\mathfrak m} \cap [\alpha-1/2H, \alpha+1/2H]} |S_{f}(\beta;X)|^2\ d\beta \ll \eta^{\frac{3}{\delta}} K^{-\frac{2-\delta}{\delta}} F^2 X.
\end{equation}
\item[(iii)] (Large values estimate)  For any $J \geq 1$, and any $1/H$-separated set of elements $\alpha_1,\dots,\alpha_J$ of $\T$, one has
\begin{equation}\label{lve-f}
\sum_{j=1}^J \int_{[\alpha_j-1/2H, \alpha_j+1/2H]} |S_{f}(\alpha;X)|^2\ d\alpha \ll K J^{1/2-\delta} F^2 X
\end{equation}
and
\begin{equation}\label{lve-g}
\sum_{j=1}^J \int_{[\alpha_j-1/2H, \alpha_j+1/2H]} |S_{g}(\alpha;X)|^2\ d\alpha \ll K J^{1/2-\delta} G^2 X.
\end{equation}
\end{itemize}
Then for all but $O_\delta(\eta H)$ integers $h$ with $|h-h_0| \leq H$, one has
\begin{equation}\label{fg2}
 \sum_n f(n) \overline{g(n+h)} = \operatorname{MT}_{h} + O( \eta FG X ).
\end{equation}
\end{lemma}

In practice, we will apply this lemma with $K = \eta^{-1} = \log^{o(1)} X$ for some slowly decaying $o(1)$, while keeping $\delta$ at a fixed value such as $1/4$; as such, the precise values of the exponents $\frac{3}{\delta}$ and $\frac{2-\delta}{\delta}$ appearing in \eqref{fF2} will not be of major significance. Also we will not apply this lemma directly to $d_k$ but rather to a truncated version for which we can establish the hypotheses (ii) and (iii) of the lemma.

\begin{proof}[Proof of Lemma~\ref{chebyshev-2}]  By translating $g$ by $h_0$, we may normalize $h_0=0$; by dividing $f,g$ by $F,G$ we may assume that $F=G=1$.  

From \eqref{sam} one has
$$
\sum_n f(n) \overline{g}(n+h) = \int_{\mathfrak M} S_{f}(\alpha;X) \overline{S_{g}(\alpha;X)} e(\alpha h)\ d\alpha + \int_{\mathfrak m} S_f(\alpha;X) \overline{S_g(\alpha;X)} e(\alpha h)\ d\alpha.$$
In view of \eqref{flip} and Chebyshev's inequality, it suffices to show that 
\begin{equation}
\label{sfg3} S \coloneqq \sum_{|h| \leq H} \left|\int_{\mathfrak m} S_{f}(\alpha;X) \overline{S_{g}(\alpha;X)} e(\alpha h)\ d\alpha\right|^2 \ll_\delta \eta^3 H X^2.
\end{equation}
By \cite[Proposition 3.1]{mrt-corr} we have an estimate for $S$:
$$ S \ll H \int_{\mathfrak m} |S_{f}(\alpha;X)| |S_{g}(\alpha;X)| \int_{\mathfrak m \cap [\alpha-1/2H, \alpha+1/2H]} |S_{f}(\beta;X)| |S_{g}(\beta;X)|\ d\beta d\alpha.$$
We partition $\T$ into $H$ intervals $I_j = [\frac{j-1}{H}, \frac{j}{H})$ of length $\frac{1}{H}$, and conclude that
\[
\begin{split} S &\ll H \sum_{i = 1}^H \left(\int_{\mathfrak m \cap I_i} |S_{f}(\alpha;X)| |S_{g}(\alpha;X)|\ d\alpha\right) \sum_{j = i-1}^{i+1} \left(\int_{\mathfrak m \cap I_j} |S_{f}(\beta;X)| |S_{g}(\beta;X)|\ d\beta\right) \\ 
&\ll H \sum_{i=1}^H \left(\int_{\mathfrak m \cap I_i} |S_{f}(\alpha;X)| |S_{g}(\alpha;X)|\ d\alpha\right)^2,
\end{split}
\]
where we used the inequality $|ab| \ll |a|^2 + |b|^2$, and extend $I_j$ to $j=0$ and $j=H+1$ in the obvious fashion. By the Cauchy-Schwarz inequality, we therefore have the bound
$$ S \ll H \sum_{i=1}^H a(I_i) b(I_i),$$
where
\begin{align*}
a(I) \coloneqq \int_{\mathfrak m \cap I} |S_{f}(\alpha;X)|^2\ d\alpha \quad \text{and} \quad b(I) \coloneqq \int_{I} |S_{g}(\alpha;X)|^2\ d\alpha.
\end{align*}
Applying H\"older's inequality, we conclude that
$$ S \ll H \|a\|_{\ell^\infty}^\delta \|a\|_{\ell^{2-\delta}}^{1-\delta} \|b\|_{\ell^{2-\delta}}.$$
From \eqref{fF2}, we have
$$ \|a\|_{\ell^\infty} \ll \eta^{\frac{3}{\delta}} K^{-\frac{2-\delta}{\delta}} X.$$
From \eqref{lve-f}, we have
$$ a(I_{i_1}) + \dots + a(I_{i_J}) \ll K J^{1/2-\delta} X $$
for any distinct intervals $I_{i_1},\dots,I_{i_J}$.  In particular, for any $1 \leq J \leq H$, the $J^{\operatorname{th}}$ largest value of $a(I_i)$ is $O( K J^{-1/2-\delta} X )$.  Since $(2-\delta) (1/2+\delta) > 1$, we thus conclude that
$$ \|a\|_{\ell^{2-\delta}} \ll_\delta K X,$$
and similarly
$$ \|b\|_{\ell^{2-\delta}} \ll_\delta K X.$$
The claim \eqref{sfg3} follows.
\end{proof}

We will need to apply Lemma~\ref{chebyshev-2} to truncated versions of $\Lambda$ and $d_k$ which we define next. First, we write $\tilde \Lambda$ for the restriction of $\Lambda$ to the primes, thus $\tilde \Lambda(n) \coloneqq \log n$ when $n$ is prime, and $\tilde \Lambda(n) = 0$ otherwise. When studying divisor functions $d_k$, in order to get good large values estimates and minor arc estimates it is convenient to truncate to a set ${\mathcal S}_{k,X}$ of ``typical'' numbers $n$.  More precisely, for any $X \geq 100, k \geq 2$ and a small fixed $\varepsilon'>0$, we let ${\mathcal S}_{k,X}$ denote the set of all $n \in (X,2X]$ satisfying the following two conditions.

\begin{itemize}
\item[(A)] Let (cf. \cite{mr})
\begin{align*}
P_1 \coloneqq (\log X)^{\psi(X)} \ & , \ Q_1 \coloneqq (\log X)^{10 k \log k} \\
P_2 \coloneqq \exp((\log\log X)^2) \ & , \ Q_2 \coloneqq \exp((\log\log X)^{5/2}) \\
P_3 \coloneqq \exp((\log X)^{3/4}) \ & , \ Q_3 \coloneqq \exp((\log X)^{5/6}).
\end{align*}
with $\psi(X) \to 0$ very slowly with $x \to \infty$.
Writing $n = a b c d$ with all the primes factors of $a,b,c$ respectively in $[P_1, Q_1]$, $[P_2, Q_2]$ and $[P_3, Q_3]$ and with $d$ having no prime factors in these intervals, we require that $a,b,c$ are greater than $1$ and square-free. 
\item[(B)]  The total number of prime factors of $n$ does not exceed $\nopf k \log\log X$, and the total number of prime factors of $n$ in the range $[X^{1/(\log\log X)^2}, 2X]$ does not exceed $10 k \log\log\log X$. 
\end{itemize}

The first condition (A) will be used to deploy a variant of the theorem \cite{mr} of the first two authors to study the behavior of $d_k(n)$ in almost all short intevals. The second condition will be used to construct an efficient point-wise majorant for $d_k(n)$ on the set of integers obeying the condition (B). Moreover the fact that the total number of prime factors of such an integer does not exceed $\nopf k \log\log X$ will be used to get a rough but useful upper bound for $d_k(n)$ at various points.

Let us start by showing that correlations of $d_k(n)$ can be well approximated by correlations of $d_k(n) 1_{{\mathcal S}_{k , X}}(n)$. 

\begin{lemma}\label{slip} Let $X \geq H \geq 100$, $k, l \geq 2$.  Let $f, \tilde f: \Z \to \R$ be the functions $f(n) \coloneqq d_k(n) 1_{(X,2X]}(n)$ and $\tilde f(n) \coloneqq d_k(n) 1_{{\mathcal S}_{k,X}}(n)$.  Then
$$
\sum_{0 < |h| \leq H} \sum_n |f(n)-\tilde f(n)| d_l(n+h) = o_{k, l, \psi}(H X \log^{k+l-2} X).
$$

If $f(n) \coloneqq \Lambda(n) 1_{(X,2X]}(n)$ and $\tilde f(n) \coloneqq \tilde \Lambda(n) 1_{(X,2X]}(n)$, we have
$$
\sum_{|h| \leq H} \sum_n |f(n)-\tilde f(n)|d_l(n+h) \ll H X^{3/4}.
$$
\end{lemma}

\begin{proof}
The second claim is trivial since the summand is always at most $X^{1/10}$ and can be non-zero only for $n \leq (2X)^{1/2}$.

To prove the first claim, let us define, for $j = 1, 2, 3$, the multiplicative functions $g_j$ by setting
\[
g_j(p^\nu) \coloneqq 
\begin{cases}
d_k(p^\nu) &\text{if $p \not \in [P_j, Q_j]$ or $\nu > 1$;} \\
0 &\text{otherwise.}
\end{cases}
\]
and, for $j = 4, 5$, the functions $g_j$ by setting $g_4(n) \coloneqq d_k(n) 1_{\Omega(n) \geq \nopf k \log \log X}$ and $g_5(n) \coloneqq d_k(n)$ if $n$ has more than $10 k \log \log \log X$ prime factors in the range $[X^{1/(\log \log X)^2}, 2X]$, and $g_5(n) = 0$ otherwise. Then it suffices by the triangle inequality to show that, for each $j = 1, \dotsc, 5$, we have
\[
\sum_{0 < |h| \leq H} \sum_{X < n \leq 2X} g_j(n)d_l(n+h) = o_{k, l, \psi}(H X \log^{k+l-2} X).
\]
For $j = 1, 2, 3$, Lemma~\ref{henriot}(iv) with $q = 1$ gives
\[
\begin{split}
\sum_{0 < |h| \leq H} \sum_{X < n \leq 2X} g_j(n)d_l(n+h) &\ll_{k,l} HX(\log X)^{k+l-2} \prod_{\substack{p \in [P_j, Q_j]}} \frac{\left(1-\frac{1}{p}\right)}{\left(1+\frac{k-1}{p}\right)} \\
&= o_{k, l, \psi}(HX(\log X)^{k+l-2}).
\end{split}
\]

When $j = 4$, we estimate, for some $0 < \delta < 1$
\[
\begin{split}
&\sum_{0 < |h| \leq H}\sum_{X < n \leq 2X} g_4(n)d_l(n+h) \\
&\leq \sum_{0 <|h| \leq H} \sum_{n \leq 2X} d_k(n) (1+\delta)^{\Omega(n) - \nopf k \log \log X} d_l(n+h)\\
& \ll_{k, l} (\log X)^{-\nopf k \log (1+\delta)} HX \prod_{p \leq 2X} \left(1+\frac{(1+\delta)k-1}{p}\right)\left(1+\frac{l-1}{p}\right) \\
&\ll HX (\log X)^{k+l-2 + [\delta - \nopf \log(1+\delta) ] k} \\
&= o(HX(\log X)^{k+l-2}).
\end{split}
\]
when $\delta$ is chosen to be a small multiple of $\varepsilon'$, and where we again applied Lemma~\ref{henriot}(iv). Finally for the case $j = 5$, let us write $\Omega_{[P, Q]}(n) = \sum_{\substack{P \leq p \leq Q \\ p^{\alpha} || n}} \alpha$. Then
\[
\begin{split}
&\sum_{0 < |h| \leq H} \sum_{X < n \leq 2X} g_5(n) d_l(n+h) \\
&\leq \sum_{0 < |h| \leq H} \sum_{n \leq 2X} d_k(n) 2^{\Omega_{[X^{1/(\log \log X)^2}, 2X]}(n) - 10 k \log \log \log X} d_l(n+h) \\
&\ll_{k, l} (\log \log X)^{-10 k \log 2} HX \prod_{p \leq X^{1/(\log \log X)^2}} \left(1+\frac{k-1}{p}\right)\\
& \quad \cdot\prod_{X^{1/(\log \log X)^2} \leq p \leq 2X} \left(1+\frac{2k-1}{p}\right) \prod_{p \leq 2X} \left(1+\frac{l-1}{p}\right)  \\
&\ll HX (\log X)^{k+l-2} (\log \log X)^{-4k},
\end{split}
\]
again by Lemma~\ref{henriot}(iv).
\end{proof}

Given parameters $Q \geq 1$ and $\delta > 0$, define the major arcs
$$ {\mathfrak M}_{Q,\delta} \coloneqq  \bigcup_{1 \leq q \leq Q} \bigcup_{a: (a,q) = 1} \left[\frac{a}{q} - \delta, \frac{a}{q} + \delta\right],$$
where we identify intervals such as $[\frac{a}{q}-\delta, \frac{a}{q}+\delta]$ with subsets of the unit circle $\T$ in the usual fashion. We will need the following major arc estimate.

\begin{proposition}[Major arc estimate]\label{major}
Let $A \geq 1$ and $k,l \geq 2$ be fixed integers, and suppose that $X \geq 2$. Write $f_k(n) \coloneqq d_k(n) 1_{{\mathcal S}_{k,X}}(n)$ and $g(n) \coloneqq \tilde\Lambda(n) 1_{(X,2X]}(n)$. Let $0 < h < X^{1-\varepsilon}$, and let $C_{k, l, h}$ and $C_{k, h}$ be as in Theorem \ref{uncond}. Then
\begin{itemize}
\item[(i)]  (Major arcs for divisor correlation conjecture)
\begin{align*}
&\int_{{\mathfrak M}_{\log^{A} X,X^{-1} \log^{3A} X}} S_{f_k}(\alpha;X) \overline{S_{f_l}(\alpha;X)} e(\alpha h)\ d\alpha = C_{k,l,h} X \log^{k+l-2} X + o_{k,l,A}(X \log^{k+l-2} X ).
\end{align*}
\item[(ii)]  (Major arcs for higher order Titchmarsh problem)
\begin{align*}
&\int_{{\mathfrak M}_{\log^{A} X,X^{-1} \log^{3A} X}} S_{g}(\alpha;X) \overline{S_{f_k}(\alpha;X)} e(\alpha h)\ d\alpha = C_{k,h} X \log^{k-1} X + o_{k,A}(X \log^{k-1} X ).
\end{align*}
\end{itemize}
\end{proposition}

We will prove this in Section~\ref{se:Major}. Notice that the corresponding claims without truncating the functions $d_k$ and $\Lambda$ (and with a more precise main term and error term) were essentially shown in \cite[Proposition 3.3]{mrt-corr}.

We will complement these major arc estimates with the following minor arc and large values estimates:

\begin{proposition}[Minor arcs for higher order divisor functions]\label{minor-arc}  Let $k \geq 2 , A \geq 1000 k \log k$, and let $H$ be such that $\log^{10 000 k \log k} X \leq H \leq X^{1-\eps}$.    Set ${\mathfrak m} \coloneqq \T \backslash {\mathfrak M}_{\log^{A} X,X^{-1} \log^{3A} X}$ and $f \coloneqq d_k(n) 1_{{\mathcal S}_{k,X}}(n)$.
Then for any $\alpha \in {\mathfrak m}$, we have
\begin{equation}\label{minor}
\int_{{\mathfrak m} \cap [\alpha-1/2H, \alpha+1/2H]} |S_{f}(\beta;X)|^2\ d\beta \ll_{k,\eps} P_1^{-1/10} X \log^{2(k-1)} X.
\end{equation}
\end{proposition}

\begin{proposition}[Large values estimates]\label{lve-prop}  Let $k \geq 2$ and $\eps>0$, let $X \geq 2$, and let $H$ be such that $\log^{10 000 k \log k} X \leq H \leq X^{1-\eps}$.  Let $\alpha_1,\dots,\alpha_J$ be a $1/H$-separated subset of $\T$.
\begin{itemize}
\item[(i)]  Let $f: \Z \to \R$ be the function $f(n) \coloneqq d_k(n) 1_{{\mathcal S}_{k,X}}(n)$. There exists a function $\psi_1(x)$, not depending on the function $\psi(x)$ used in the definition of $P_1$, such that $\psi_1(x) \to 0$ with $x\to \infty$ and
\begin{equation}\label{fx}
\sum_{j=1}^J \int_{[\alpha_j-1/2H,\alpha_j+1/2H]} |S_f(\beta;X)|^2\ d\beta \ll_{k,\eps} J^{1/4} X \log^{2(k-1)+\psi_1(x)} X.
\end{equation}
\item[(ii)]  If instead $f: \Z \to \R$ is the function $f(n) \coloneqq \tilde \Lambda(n) 1_{(X,2X]}(n)$, then one has
$$
\sum_{j=1}^J \int_{[\alpha_j-1/2H,\alpha_j+1/2H]} |S_f(\beta;X)|^2\ d\beta \ll_{k,\eps} J^{1/4} X \log^{o(1)} X.
$$
\end{itemize}
\end{proposition}

Proposition \ref{minor-arc} will be proven in Section \ref{sec:minor}.  The exponent of $-1/10$ appearing in \eqref{minor} is not optimal, but any negative exponent would suffice for our argument here. We prove Proposition \ref{lve-prop} in Section \ref{lve-sec}.  The exponent of $1/4$ on the right-hand side is not important; any exponent between $0$ and $1/2$ would have sufficed.

\begin{proof}[Proof of Theorem~\ref{uncond} assuming Propositions \ref{major}--\ref{lve-prop}]
Let us begin with part (i). Let $\eps>0$, $k \geq l \geq 2$, $X \geq 2$, and let $H$ be such that $\log^{10000 k \log k} X \leq H \leq X^{1-\eps}$. Let $\psi_1(x)$ be as in the previous proposition and let $\psi(X) \geq 1000 \psi_1(X)$ be tending to $0$ very slowly with $X \to \infty$. We apply Lemma \ref{chebyshev-2} with 
\begin{align*}
f(n) &\coloneqq d_k(n) 1_{{\mathcal S}_{k,X}}(n), \quad g(n) \coloneqq d_l(n) 1_{{\mathcal S}_{l,X}}(n), \\
{\mathfrak M} &\coloneqq {\mathfrak M}_{\log^{1000k \log k} X,X^{-1} \log^{3000 k \log k} X}, \qquad
{\mathfrak m} \coloneqq \T \backslash {\mathfrak M}, \\
F &\coloneqq \log^{k-1} X, \qquad G \coloneqq \log^{l-1} X, \\
\delta &\coloneqq 1/4, \qquad \eta \coloneqq P_1^{-1/500}, \qquad K \coloneqq P_1^{1/500}, \\
\operatorname{MT}_h &\coloneqq \int_{\mathfrak M} S_{f}(\alpha;X) \overline{S_{g}(\alpha;X)} e(\alpha h)\ d\alpha,\\
h_0 &\coloneqq 0.
\end{align*}
The error between $\sum_n f(n) g(n+h)$ and $\sum_{X < n \leq 2X} d_k(n) d_l(n+h)$ is acceptable for almost all $|h| \leq H$ by Lemma~\ref{slip}. Furthermore the error between $C_{k, l, h} X \log^{k+l-2} X$ and $\operatorname{MT}_h$ is acceptable for all $h \neq 0$ by Proposition~\ref{major}. 

Hence it suffices to verify the hypotheses (i)--(iii) of Lemma \ref{chebyshev-2}. The major arc estimate \eqref{flip} is trivial. The minor arc estimate \eqref{fF2} follows from Proposition \ref{minor-arc}, and the large values estimates \eqref{lve-f}, \eqref{lve-g} follow from Proposition \ref{lve-prop}. This proves Theorem \ref{uncond}(i).

The proof of part (ii) is similar, replacing $l$ by $1$, $d_l 1_{{\mathcal S}_{l,X}}$ by $\tilde \Lambda 1_{(x, 2X]}$, and $\operatorname{MT}_h$ by the corresponding major arc integral; we leave the details to the interested reader.
\end{proof}

\section{Major arc estimates}
\label{se:Major}

In this section we prove Proposition~\ref{major}. 
Before proving the proposition we start with a somewhat general lemma related to the circle method. Let us first define a class of functions to which the lemma applies.

\begin{definition}
\label{def:Cclass}
For $Q, Y \in [2, X]$, $F \geq 1$ and $\widetilde{f} \colon \mathbb{N} \to \mathbb{R}_{\geq 0}$ with $\widetilde{f}(d) \ll d^{A}$ for some fixed constant $A > 0$, let $\mathcal{C}(X, Y, Q, \widetilde{f}, F)$ denote the class of functions $f \colon \mathbb{N} \to \mathbb{C}$ such that the following two conditions hold.
\begin{itemize}
\item[(i)] For all $d \leq Q$, all $Y' \in [Y, X]$ and all $x \in [\frac{X}{2d}, \frac{4X}{d}]$, one has
$$
\frac{1}{Y'} \sum_{x \leq n \leq x + Y'} f(d n) \leq \widetilde{f}(d) \cdot F.
$$
\item[(ii)] For all non-principal characters of modulus $q \leq Q$, all $Y' \in [Y, X]$, all $d \leq Q$, and $x \in [\frac{X}{2d}, \frac{4X}{d} ]$, one has 
\[
\frac{1}{Y'} \sum_{x \leq n \leq x + Y'} f(d n) \chi(n) \ll_{K} F Q^{-K}
\]
for any $K \geq 1$.
\end{itemize}
\end{definition}

\begin{lemma} \label{softarclemma}
For $X \geq Q, Y \geq 2$, $F, G \geq 1$ and $\widetilde{f}, \widetilde{g} \colon \mathbb{N} \to \mathbb{C}$ with $\widetilde{f}(d), \widetilde{g}(d) \ll d^{A}$ for some fixed constant $A > 0$.
Let $f,g : \mathbb{N} \rightarrow \mathbb{R}^{+}$ be supported on $(X, 2X]$, and suppose that $f \in \mathcal{C}(X, Y/2, Q, \widetilde{f}, F)$ and $g \in \mathcal{C}(X, Y/2, Q, \widetilde{g}, G)$. 
Then, for any $h$, 
\begin{align*}
\int_{\mathfrak{M}_{Q,1/Y}} & S_{f}(\alpha;X) \overline{S_g(\alpha;X)} e(h \alpha) d \alpha \\ & \ll 
X FG \sum_{q \leq Q} |c_q(h)|  \Big ( \sum_{d | q} \frac{\mu(q / d)^2 \widetilde{f}(d)}{\varphi(q/d) d}\Big ) \Big ( \sum_{d | q} \frac{\mu(q / d)^2 \widetilde{g}(d)}{\varphi(q/d) d} \Big )+ O_K(FGXQ^{-K})
\end{align*}
for any $K \geq 1$, where
\[
c_q(h) \coloneqq \sum_{\substack{1 \leq b \leq q \\ (b, q) = 1}} e\left(\frac{hb}{q}\right)
\]
is the Ramanujan sum.
\end{lemma}
\begin{proof}
Consider
$$
S_{f} \Big ( \frac{a}{q} + \beta;X \Big ) = \sum_{X < n \leq 2X} f(n) e \Big ( \frac{n a}{q} \Big ) e(n \beta)
$$
with $q \leq Q$, $(a,q) = 1$, and $|\beta| \leq 1/Y$. Splitting the summation according to $(n, q)$ and $n/(n,q) \pmod{q/(n, q)}$, we see that
$$
S_{f} \Big ( \frac{a}{q} + \beta;X \Big ) = \sum_{d | q} \sum_{\substack{\ell \pmod{q/d} \\ (\ell, q / d) = 1}} e \Big ( \frac{\ell a}{q / d} \Big ) 
\sum_{\substack{X/d < n \leq 2X / d \\ n \equiv \ell \pmod{q / d}}} f(d n) e (d n \beta).
$$
Expressing the congruence condition in terms of Dirichlet characters and separating the contribution of the principal character, we obtain
\[
\begin{split}
S_{f} \Big ( \frac{a}{q} + \beta;X \Big ) &=  
\sum_{d | q}  \sum_{\substack{\ell \pmod{q/d} \\ (\ell, q / d) = 1}} e \Big ( \frac{\ell a }{q / d} \Big ) 
\frac{1}{\varphi(q / d)} \sum_{\substack{X/d < n \leq 2X/d \\ (n, q/d) = 1}} f(d n) e( d n \beta) \\
& \quad + \sum_{d | q}  \sum_{\substack{\ell \pmod{q/d} \\ (\ell, q / d) = 1}} e \Big ( \frac{\ell a}{q / d} \Big ) \frac{1}{\varphi(q/d)} \sum_{\chi \neq \chi_{0} \pmod{q / d}} \widetilde{\chi}(\ell) \sum_{X/d < n \leq 2X/d} f( d n) \chi(n) e( d n \beta)  \\
&=: M_f(a, q; \beta) + E_f(a, q; \beta),
\end{split}
\]
say (see also e.g. \cite[Lemma 3.2]{breteche1} or \cite[Lemma 2.2]{breteche2}). Notice that, since $(a,q) = 1$, we have 
$$
\sum_{\substack{\ell \pmod{q/d} \\ (\ell, q / d) = 1}} e \Big ( \frac{\ell a}{q / d} \Big ) = c_{q/d}(a) = \mu(q/d). 
$$
Therefore, 
\begin{align*}
M_{f}(a, q; \beta) & = \sum_{d | q} \frac{\mu(q/d)}{\varphi(q/d)} \sum_{\substack{X/d < n \leq 2X/d \\ (n,q/d) = 1}} f(d n) e(d n \beta) =  \sum_{d | q} \frac{\mu(q/d)}{\varphi(q/d)} \widetilde{M}_f(d; \beta),
\end{align*}
say. Hence
\begin{equation}
\label{eq:SftildeME}
\begin{split}
\int_{\mathfrak{M}(Q,1/Y)} S_{f}(\alpha;X) \overline{S_g(\alpha;X)} e(h \alpha) d \alpha &= \sum_{q \leq Q} \sum_{(a, q) = 1} \int_{|\beta| \leq 1/Y} \left(\sum_{d \mid q} \frac{\mu(q/d)}{\varphi(q/d)} \widetilde{M}_f(d; \beta) + E_f(a, q; \beta)\right) \\
& \quad \cdot \left(\sum_{d \mid q} \frac{\mu(q/d)}{\varphi(q/d)} \widetilde{M}_g(d; \beta) + E_g(a, q; \beta)\right) e\left(h\left(\frac{a}{q}+\beta\right)\right) d\beta
\end{split}
\end{equation}
Applying Gallagher's lemma (Lemma \ref{le:Gallagher}) and Definition \ref{def:Cclass}(i) for $f$, we see that 
\[
\int_{|\beta| \leq 1/Y} |\widetilde{M}_f(d; \beta)|^2 d \beta \ll \frac{\widetilde{f}(d)^2 F^2 X}{d^2},
\]
Furthermore
\[
\begin{split}
&\sum_{q \leq Q} \sum_{(a,q) = 1} \int_{|\beta| \leq 1/Y} |E_f(a, q; \beta)|^2 d \beta \\
&\ll Q^{10} \cdot \sup_{\substack{ \chi \neq \chi_0 \pmod{q} \\ d, q \leq Q}} \int_{|\beta| \leq 1/Y} \Big | \sum_{X/d < n \leq 2X/d} f(d n) \chi(n)e(d n \beta) \Big |^2 d \beta \\
&\ll_{K} X F^2 \cdot Q^{-K}
\end{split}
\]
for any $K \geq 1$, where we have applied Gallagher's lemma (Lemma \ref{le:Gallagher}) and Definition \ref{def:Cclass}(ii) for $f$. Similar bounds hold with $f, \widetilde{f}$ and $F$ replaced by $g, \widetilde{g}$ and $G$.

The claim follows now by multiplying \eqref{eq:SftildeME} out, using Cauchy-Schwarz, and then applying the previous bounds.
\end{proof}

For the proof of Proposition \ref{major} we will require the following major arc estimate. 
\begin{lemma} \label{majorarcmain} Let $k \geq \ell \geq 2$ and $h \in [1, X^{1-\varepsilon}]$ be integers. Then
\begin{equation*}
\int_{\mathfrak{M}_{\log^{A} X , X^{-1} \log^{3A} X}} S_{d_k}(\alpha;X) \overline{S_{d_{\ell}}(\alpha;X)} e(h \alpha) d \alpha = X P_{k,\ell,h}(\log X) + O_{A}(X (\log X)^{-A/3})
\end{equation*}
for any $A > 0$, 
with $P_{k,\ell,h}$ a polynomial of degree $k + \ell - 2$. 

We also have,
$$
\int_{\mathfrak{M}_{\log^{A} X, X^{-1} \log^{3 A} X}} S_{\Lambda}(\alpha;X) \overline{S_{d_{\ell} }(\alpha;X)} e(h \alpha) d\alpha = X Q_{\ell,h}(\log X) + O_{A}(X (\log X)^{-A / 3})
$$
for any $A > 0$, 
with $Q_{\ell,h}$ a polynomial of degree $\ell - 1$. 
\end{lemma}
\begin{proof}
See \cite[Proposition 3.3 ii) and iii)]{mrt-corr}.  
\end{proof}

Now we are ready to prove Proposition~\ref{major}. Let us start with the case of divisor correlations. Given Lemma \ref{majorarcmain} it suffices to show that
$$
\int_{\mathfrak{M}_{\log^{A} X, X^{-1} \log^{3A} X}} \Big ( S_{f_k}(\alpha;X) \overline{S_{f_{\ell}}(\alpha;X)} - S_{d_k}(\alpha;X) \overline{S_{d_{\ell}}(\alpha;X)} \Big ) e(h \alpha) d \alpha = o(X (\log X)^{2k - 2}).
$$

Now since 
$$
S_{f_k}(\alpha;X) = S_{d_k}(\alpha;X) - S_{g_k}(\alpha;X)
$$
with $g_k(n) \coloneqq d_k(n) \mathbf{1}_{n \in (X, 2X] \setminus \mathcal{S}_{k, X}}$, it suffices to show that
$$
\int_{\mathfrak{M}_{\log^{A} X , X^{-1} \log^{3A} X}} S_{f}(\alpha;X) \overline{S_{g}(\alpha;X)} e(h \alpha) d \alpha  = o (X (\log X)^{2k - 2})
$$
with $f = f_k$ and $g = g_l$, with $f = g_k$ and $g = f_l$ and with $f = g_k$ and $g = g_l$. Let us concentrate on the first possibility.

We apply the previous lemma with $Q = (\log X)^{A}$ and $Y = X / (\log X)^{3A}$. Furthermore we pick $\widetilde{f}(n) = d_k(n)$ and $\widetilde{g}(n) = d_l(n)$. A simpler variant of the proof of Lemma~\ref{slip} (and the inequality $d_k(ef) \leq d_k(e) d_k(f)$) shows that Definition~\ref{def:Cclass}(i) is satisfied with
$$
F \ll (\log X)^{k - 1} \text{ and } G = o((\log X)^{k - 1}).
$$
Thus to conclude it remains to check that Definition~\ref{def:Cclass}(ii) holds for our choice of $f$ and $g$, and we prove that this is indeed the case in Lemma \ref{lemmaCheck} below.

In the case of correlations of the divisor function with the von Mangoldt function we notice that, by Lemma \ref{majorarcmain},
$$
\int_{\mathfrak{M}_{\log^{A} X , X^{-1} \log^{3A} X}} S_{\Lambda }(\alpha;X) \overline{S_{d_{\ell}}(\alpha;X)} e(h \alpha) d \alpha = X Q_{\ell,h}(\log X) + O_{A}(X (\log X)^{-A/3})
$$
and by the same argument as above it is therefore enough to show that, 
$$
\int_{\mathfrak{M}_{\log^{A} X , X^{-1} \log^{3A} X}} S_{f}(\alpha;X) \overline{S_{g_{\ell}}(\alpha;X)} e(h \alpha) d \alpha  = o (X (\log X)^{2k - 2})
$$
with $f = \widetilde{\Lambda} 1_{(X, 2X]}$ and $g_{\ell} = d_{\ell} 1_{(X, 2X] \backslash \mathcal{S}_{\ell, X}}$ and other similar cases. We then apply Lemma \ref{softarclemma} with $Q  = (\log X)^{A}$ and $Y = X / (\log X)^{3 A}$. We pick $\widetilde{g}(n) = d_{\ell}(n)$ and $\widetilde{f}(n) = 1_{n = 1}$. Standard sieve bounds for primes give that Definition~\ref{def:Cclass}(i) holds with $F \ll 1$ and a variant of Lemma \ref{slip} gives that it holds with $G = o((\log X)^{\ell - 1})$. Thus again it remains to check that Definition~\ref{def:Cclass}(ii) holds for $f$ and $g$. The assumption for $f$ follows from Siegel-Walfisz. The assumption for $g$ is verified in the following lemma. 

\begin{lemma} \label{lemmaCheck}
Let $A, B, K \geq 1$. Let $d \leq (\log X)^{A}$ and let $\chi$ be a non-principal character of conductor at most $(\log X)^{B}$. Then, uniformly for any $\alpha \in [1, 4]$,
$$
\sum_{\substack{d n \in \mathcal{S}_{k, X} \\ X/d < n \leq \alpha X/d}} d_k(d n) \chi(n) \ll_{A, B, K} X (\log X)^{-K}
$$
for any $K \geq 1$.
The same bound holds if the condition $d n \in \mathcal{S}_{k,X}$ is replaced with $d n \in (X, 2X] \setminus \mathcal{S}_{k, X}$. 
\end{lemma}

\begin{proof}
  First let us focus on the case with $d n \in \mathcal{S}_{k, X}$. Membership in $\mathcal{S}_{k, X}$ amounts to requiring that the integer $m$ in question is such that $\Omega(m; I_{j}) \geq 1$ for $j = 1,2,3$ where $I_j = [P_j, Q_j]$, that $m$ has no repeated prime factors in intervals $I_{j}$, and finally that $\Omega(m) \leq (1 + \varepsilon') k \log\log X$ and $\Omega(m; J) \leq 10 k \log\log\log X$ where $J = [X^{1/(\log\log X)^2}, 2X]$.

  Therefore,
  \begin{align*}
    \sum_{\substack{d n \in \mathcal{S}_{k, X} \\X / d < n \leq 2 X / d}} d_k(d n) \chi(n) =  
     \sum_{\substack{1 \leq \ell, \ell_1, \ell_2, \ell_3 \leq (1 + \varepsilon') k \log\log X \\ \ell_4 \leq 10 k \log\log\log X}} \sum_{\substack{\Omega(d n) = \ell \\ \Omega(d n; I_j) = \ell_j \ , \ j = 1,2,3 \\ \Omega(d n; J) = \ell_4}} d_k(d n) \chi(n) \prod_{j = 1}^{3} \mu^2(d n; I_j)
  \end{align*}
where for any interval $J$, the function $\mu^2(n; J)$ is defined by setting $\mu^2(n; J) \coloneqq 1$ if all the prime factors of $n$ lying in $J$ have multiplicity $1$, and $\mu^2(n; J) \coloneqq 0$ otherwise.
  
  Using Cauchy's formula we can express the preceding expression as
  $$
  \sum_{\substack{1 \leq \ell, \ell_1, \ell_2, \ell_3 \leq (1 + \varepsilon')k \log\log X \\ \ell_4 \leq 10 k \log\log\log X}} \frac{1}{(2\pi i)^{5}} \oint \ldots \oint S_{k, d, \chi}(X; z, z_1, z_2, z_3, z_4) \cdot \frac{dz dz_1 dz_2 dz_3 dz_4}{z^{\ell + 1} z_1^{\ell_1 + 1} z_2^{\ell_2 + 1} z_3^{\ell_3 + 1} z_4^{\ell_4 + 1}} 
  $$
  where we integrate over $|z| = |z_1| = |z_2| = |z_3| = |z_4| = 1$ and 
  where
  $$
  S_{k, d, \chi}(X; z, z_1, z_2, z_3, z_4) \coloneqq \sum_{\substack{X / d < n \leq \alpha X / d}} d_k(d n) \chi(n) z^{\Omega(d n)} \Big ( \prod_{j = 1}^{3} z_j^{\Omega(d n; I_j)} \mu^2(dn ; I_j) \Big )  z_4^{\Omega(d n; J)}
  $$
  Strictly speaking our notation for $S_{k,d}(X; \cdot)$ should also involve the intervals $I_j, J$, but to avoid extra clutter this dependence will not be made explicit in the notation. 
  To conclude the argument, it is enough to show that uniformly in $|z| = |z_j| = 1$ ($j = 1,2,3,4$) and $d \leq (\log X)^{A}$ and non-principal characters $\chi$ of conductor $\leq (\log X)^{B}$ we have
  $$
  S_{k, d, \chi}(X; z, z_1, z_2, z_3, z_4) \ll_{A} X (\log X)^{-A}. 
  $$
  To accomplish this we notice that, denoting the coefficients of $S_{k, d, \chi}(X; z, z_1, z_2, z_3, z_4)$ by $d_k(d n) z^{\Omega(d n)} a_{k, d}(n; z_1, z_2, z_3, z_4) \chi(n)$, we have,
  $$
  S_{k, d, \chi}(X; z, z_1, z_2, z_3, z_3) = \frac{1}{2\pi i} \int_{1 + \varepsilon - i X^{1 + \varepsilon}}^{1 + \varepsilon + X^{1 + \varepsilon}} L_{k, d}(s; \chi, z, z_1, z_2, z_3, z_4) \cdot 
  \frac{(\alpha^s - 1) (X / d)^s}{s} ds  + O_{\varepsilon}(X^{-\varepsilon}). 
  $$
  where
  $$
L_{k, d}(s; \chi, z, z_1, z_2, z_3, z_4) \coloneqq \sum_{n = 1}^{\infty} \frac{d_k(d n) \chi(n) z^{\Omega(d n)} a_{k, d}(n; z_1, z_2, z_3, z_4)}{n^s}. 
$$
The function involved are multiplicative at $(n, d) = 1$. Thus we see that the above Dirichlet series factors as,
$$
(B_{k, d} \cdot A_{k,d})(s, z, z_1, z_2, z_3, z_4) \cdot L(s, \chi)^{k z}
\prod_{j = 1}^{3} \prod_{p \in I_j} \Big ( 1 + \frac{k z (z_j - 1) \chi(p)}{p^s} \Big ) \prod_{p \in J} \Big ( 1 + \frac{k z (z_4 - 1) \chi(p)}{p^s} \Big ). 
$$
where
$B_{k,d}$ is absolutely convergent and $\ll_{\varepsilon} 1$ in the region $\Re s > \tfrac 12 + \varepsilon$, while $A_{k, d}$ involves only primes $p$ dividing $d$ and satisfies the bound $A_{k, d}(\cdot) \ll d^{\varepsilon}$ in the region $\Re s > 0$, uniformly in $z, z_i$ with $|z| = |z_i| = 1$.

We now shift the contour of integration to the region $\sigma = 1 - \frac{C \log\log X}{\log X}$ for some arbitrarily large but fixed $C > 0$. Since $|\log L(s, \chi)| \ll \log\log X$ in the region $\Re s > \sigma$ and $|t| \leq X^{1 + \varepsilon}$ we have $| L(s, \chi)^{k z} |  \ll (\log X)^{A}$ for some large power of $A$ depending on $k$. Moreover in this region, we have,
$$
\prod_{j = 1}^{3} \prod_{p \in I_j} \Big ( 1 + \frac{k z (z_j - 1) \chi(p)}{p^s} \Big ) \ll (\log X)^{A}
$$
for some large $A$ depending on $k$. This bound follows from the trivial bound
$$
\sum_{p \leq X^{1/(\log\log X)^2}} \frac{1}{p^{\sigma}} \ll \log\log X. 
$$
Finally, from Perron's formula we have
$$
\sum_{p \in J} \frac{\chi(p)}{p^{\sigma + it}} = \frac{1}{2\pi i} \int_{\varepsilon - i X^{1 + \varepsilon}}^{\varepsilon + i X^{1 + \varepsilon}} \log L(s + \sigma + it, \chi) \cdot ((2X)^s - X^{s / (\log\log X)^2}) \cdot \frac{ds}{s} + O(X^{-\varepsilon/2})
$$
and using the Vinogradov-Korobov zero-free region we see that this is $O_{C}(\log^{-C+2} X)$. Therefore
$$
\prod_{p \in J} \Big ( 1 + \frac{\chi(p) k z (z_4 - 1)}{p^s} \Big ) \ll 1
$$
in the region $\Re s > \sigma = 1 - \frac{C \log\log X}{\log X}$. Inserting all this information into the contour representation of $S_{k,d}(X; z, z_1, z_2, z_3, z_4)$ we conclude that
$$
S_{k, d, \chi}(X, z, z_1, z_2, z_3, z_4) \ll_{K} (X / d) \cdot d^{\varepsilon} \log^A X \cdot \log^{-C} X
$$
with $A$ depending on $k$ and $C$ arbitrarily large. Since moreover $d \leq \log^A X$ for some large fixed $A$, we conclude that $S_{k, d, \chi}(X; z, z_1, z_2, z_3, z_4) \ll_{K} X \log^{-K} X$ for any $K > 0$ as claimed.

The proof for the case $d n \not \in \mathcal{S}_{k, X}$ is essentially identical since after inclusion-exclusion and the use of Cauchy's formula the problem also boils down to showing that $S_{k,d, \chi}(X; z, z_1, z_2, z_3, z_4) \ll_{A} X \log^{-A} X$. 

\end{proof}
 
\section{Minor arc estimates}
\label{sec:minor}
The proof of Proposition~\ref{minor-arc} is reduced to estimation of mean values of Dirichlet polynomials by the following result.
\begin{lemma}\label{spem}  Let $X \geq 1$, and let $f \colon \N \to \C$ be a function supported on $(X,2X]$.  Let $q \geq 1$, let $a$ be coprime to $q$, and let $\beta, \eta$ be real numbers with $|\beta| \ll \eta \ll 1$.  Let 
\begin{equation}\label{I-def}
 I \coloneqq \left\{ t \in \R: \eta |\beta| X \leq |t| \leq \frac{|\beta| X}{\eta} \right\}
\end{equation}
Then
\begin{align*}
\int_{\beta < |\theta| \leq 2\beta} \left|S_f\left(\frac{a}{q} + \theta;X\right)\right|^2\ d\theta  &\ll \frac{d_2(q)^4}{q} 
\sup_{q=q_0 q_1}  \int_I \left( \sum_{\chi\ (q_1)} \left|\sum_n \frac{f(nq_0) \chi(n)}{n^{1/2+it}} \right|\right)^2\ dt  \\
&\quad + \left(\eta + \frac{1}{|\beta| X}\right)^2 \int_\R \left(\beta \sum_{x \leq n \leq x+1/\beta} |f(n)|\right)^2\ dx.
\end{align*}
\end{lemma}
\begin{proof}
This follows from \cite[Corollary 5.3]{mrt-corr}.
\end{proof}

Write $Q \coloneqq (\log X)^{A}$ with $A \geq 1000 k \log k$. For $\alpha \in {\mathfrak m} = \T \backslash {\mathfrak M}_{Q,Q^3/X}$, we write $\alpha = a/q + \beta$ with $q \leq Q$ and $|\beta| \leq 1/(qQ)$. Note that for any $\theta \in \mathfrak{m}$, we have $|\theta - a/q| \geq Q^3/X$. Hence
\[
\int_{{\mathfrak m} \cap [\alpha-1/2H, \alpha+1/2H]} |S_{f}(\beta;X)|^2\ d\beta \leq 2 \int_{Q^3/X}^{\frac{1}{qQ} + \frac{1}{2H}} \left|S_{f}\left(\frac{a}{q} + \beta;X\right)\right|^2\ d\beta. 
\]
We split the integration into dyadic intervals and apply the previous lemma with $\eta \coloneqq (\log X)^{-10}$. The term on the second line can be estimated squaring out, and using Lemma~\ref{henriot} for the off-diagonal terms and the upper bound $f(n) \leq (\log X)^{(1+\varepsilon') k \log k}$ for the diagonal terms. Hence it suffices to show that, when $q \leq Q, q_0 \mid q$ and $q_1 = q/q_0$, we have
\begin{equation}
\label{eq:minorArcClaim}
\int_{Q^3/(\log X)^{10}}^{X (\log X)^{10}(\frac{1}{qQ} + \frac{1}{H})} \left(\sum_{\chi \ (q_1)} \left|\sum_n \frac{f(nq_0) \chi(n)}{n^{1/2+it}} \right|\right)^2 dt \ll q P_1^{-1/8} X (\log X)^{2(k-1)}.
\end{equation}

To deal with mean values of Dirichlet polynomials without losing any log-factors, we shall use the following variant of the mean value theorem for Dirichlet polynomials.

\begin{lemma}[Log-free mean value theorem]\label{mvt}
Fix $\delta \in (0, 1)$, and let $A, B \geq 1$. Let $a(n)$ be a sequence with $|a(n)| \leq f(n)$ for some $f \in \mathcal{M}(A, B, \delta^3/1000)$. Let $X \geq Y \geq X^\delta \geq 2$. Then
\[
\begin{split}
&\sum_{\chi \pmod{q}} \int_{-T}^{T} \Big | \sum_{X < n \leq X+Y} \frac{a(n) \chi(n)}{n^{1/2 + it}} \Big |^2 dt \\
&\ll \frac{\varphi(q) T}{X} \sum_{X < n \leq X+Y} |a(n)|^2 + Y \left(\frac{\varphi(q)}{q}\right)^2 \prod_{\substack{p \leq 2X \\ p \nmid q}} \Big (1 + \frac{2 |f(p)| - 2}{p} \Big ). 
\end{split}
\]
\end{lemma}
\begin{proof}
Write $I$ for the left hand side of the claim. Let $\Phi \geq 0$ be a smooth function with $\Phi(x) \geq 1$ for $|x| \leq 1$ and with $\widehat{\Phi}(x) = 0$ for $|x| > 1$. Then
\[
\begin{split}
I & \leq \sum_{\chi \pmod{q}} \int_{\mathbb{R}} \Big | \sum_{X < n \leq X+Y} \frac{a(n) \chi(n)}{n^{1/2 + it}} \Big |^2 \cdot \Phi \Big ( \frac{t}{T} \Big ) dt \\
&= \sum_{\chi \pmod{q}} \sum_{X < m,n \leq X+Y} \frac{a(m) \overline{a(n)}}{(mn)^{1/2}} \cdot \chi(m) \overline{\chi(n)} T \widehat{\Phi} \Big ( T \log \frac{m}{n} \Big ) 
\end{split}
\]
The compact support of $\widehat{\Phi}$ ensures that upon writing $m = n + h$ only the terms with $|h| \leq 2 X / T$ survive. Moreover executing the sum over characters we find that, 
$$
\sum_{\chi \pmod{q}} \chi(n) \overline{\chi(m)} = 1_{(n m, q) = 1} \cdot 1_{n \equiv m \pmod{q}} \cdot \varphi(q). 
$$
Hence
\[
\begin{split}
I &\leq \varphi(q) \frac{T}{X} \sum_{\substack{|h| \leq 2 X / T \\ q | h}} \sum_{\substack{X < n \leq X+Y \\ (n (n + h), q) = 1}} |a(n) a(n + h)|  \\
&\leq \varphi(q) \cdot \frac{T}{X} \sum_{\substack{X < n \leq X+Y \\ (n, q) = 1}} |a(n)|^2 + 
\varphi(q) \cdot \frac{T}{X} \sum_{0 < |h| \leq 2 X / (T q)} \sum_{\substack{X < n \leq X+Y \\ (n, q) = 1}} |f(n) f(n + h q)|
\end{split}
\]
and the claim follows from Lemma~\ref{henriot}(iv).
\end{proof}

By Cauchy-Schwarz and the previous lemma, (\ref{eq:minorArcClaim}) holds trivially unless $q_0 \leq P_1^{1/2}$ which we assume from now on. Let ${\mathcal T}$ denote the interval
\[
\mathcal{T} = \mathcal{T}_{q,Q} := \left[\frac{Q^3}{(\log X)^{10}}, X (\log X)^{10}\left(\frac{1}{qQ} + \frac{1}{H}\right)\right]. 
\]
Hence, after an application of Cauchy-Schwarz, our task is to show that, for each $q \leq Q, q_0 \mid q$ and $q_1 = q/q_0$ with $q_0 \leq P_1^{1/2}$, we have
\[
\sum_{\chi \ (q_1)} \int_{\mathcal{T}} \left|\sum_n \frac{f(nq_0) \chi(n)}{n^{1/2+it}} \right|^2 dt \ll P_1^{-1/8}  X(\log X)^{2(k-1)}.
\]
To show this we use the method in~\cite{mr}.

Set $D := P_1^{1/6}$. Each Dirichlet polynomial
$$
P_j (t, \chi) := \sum_{P_j \leq p \leq Q_j} \frac{k \cdot \chi(p)}{p^{1/2 + it}} 
$$
can be decomposed as
$$
P_j(t, \chi) = \sum_{\lfloor D \log P_j \rfloor \leq v \leq D \log Q_j} Q_{v, j}(t) \ , \ Q_{v,j}(t, \chi) := \sum_{\substack{P_j \leq p \leq Q_j \\ e^{v / D} \leq p < e^{(v + 1)/D}}} \frac{k \chi(p)}{p^{1/2 + it}}. 
$$
Take $\delta_1 \coloneqq 1/4-1/100$ and $\delta_2 \coloneqq 1/4-1/50$ and, for each character $\chi$ of period $q_1$, define
\[
\begin{split}
\mathcal{T}_{1}(\chi) &= \{ t \in \mathcal{T} \colon |Q_{v,1}(t, \chi)| \leq e^{(1/2-\delta_1) v / D} \text{ for all $\lfloor D \log P_1 \rfloor \leq v \leq D \log Q_1$}\}, \\
\mathcal{T}_2(\chi) &= \{t \in \mathcal{T} \setminus \mathcal{T}_1(\chi) \colon |Q_{v,2}(t, \chi)| \leq e^{(1/2-\delta_2) v / D} \text{ for all $\lfloor D \log P_2 \rfloor \leq v \leq D \log Q_2$}\} \\
\mathcal{T}_3(\chi) &= \mathcal{T} \setminus (\mathcal{T}_1(\chi) \cup \mathcal{T}_2(\chi)).
\end{split}
\] 

Similarly to~\cite[Lemma 12]{mr}, using Lemma \ref{mvt} and Cauchy-Schwarz, we have, for $j =1, 2, 3$,
\begin{equation}
\label{eq:Ramare}
\begin{split}
& \sum_{\chi \pmod{q_1}} \int_{\mathcal{T}_j(\chi)}  \Big | \sum_{\substack{n \sim X/q_0}} \frac{f(nq_0) \chi(n)}{n^{1/2 + it}} \Big |^2 dt
\\ & \ll D \log Q_j  \sum_{\lfloor D \log P_j \rfloor \leq v \leq D \log Q_j} \sum_{\chi \pmod{q_1}} \int_{\mathcal{T}_j(\chi)} \Big | Q_{v, j}(t, \chi) R_{v, j}(t , \chi) \Big |^2 dt + X P_1^{-1/6} (\log X)^{2(k - 1)}
\end{split}
\end{equation}
where 
$$
R_{v,j}(t, \chi) := \sum_{\substack{(X/q_0) e^{-v/D} \leq m \leq 2 (X/q_0) e^{-v/D}\\ m \in \mathcal{S}^{\star}}} \frac{f(mq_0) \chi(m)}{m^{1/2 + it}} \cdot \frac{1}{\#\{P_j \leq q \leq Q_j : q | m\} + 1}
$$
where $\mathcal{S}^{\star}$ is a certain subset of $\mathcal{S}$ (which differs from ${\mathcal S}$ only through the to number of prime factors having decreased by one, and the possibility of there not being prime factors in the range $[P_j, Q_j]$). 

For $j =1, 2, 3$, we write $I_j$ for the first term on the right hand side of~\eqref{eq:Ramare}. We will treat each expression $I_1,I_2,I_3$ differently.

Using the definition of $\mathcal{T}_1(\chi)$ and the mean-value theorem (Lemma~\ref{mvt}) we bound 
\[
\begin{split}
I_1 &\ll D \log Q_1 \sum_{\lfloor D \log P_1 \rfloor \leq v \leq D \log Q_1} e^{2 (1/2- \delta_1) v / D} \\ & \times\Big ( \varphi(q_1) \left(\frac{X}{q Q} + \frac{X}{H}\right) (\log X)^{(1+\varepsilon') k \log k + k+10} + \frac{X}{e^{v/D}} (\log X)^{2k - 2} \Big ) \\
& \ll D^2 (\log Q_1)^2 P_1^{- 2 \delta_1} X (\log X)^{2 k - 2} \ll X (\log X)^{2k - 2} P_1^{-1/7},
\end{split}
\]
where we have recall that $D = P_1^{1/6}$, $H, Q > (\log X)^{1000 k\log k}$, $Q_1 = (\log X)^{10k\log k}$ and $\delta_1 = \tfrac 14-1/100$. 

Now consider $I_2$. By the definition of $\mathcal{T}_2(\chi)$
$$
I_2 \ll D \log Q_2 \sum_{\lfloor D \log P_2 \rfloor \leq v \leq D \log Q_2} e^{2 (1/2- \delta_2) v / D} \sum_{\chi \pmod{q_1}} \int_{\mathcal{T}_2(\chi)} |R_{v,2}(t, \chi)|^2 dt
$$
For each $t \in \mathcal{T}_2(\chi)$ there is a $j$ with $\lfloor D \log P_1 \rfloor \leq j \leq D \log Q_1$ such that $|Q_{j, 1}(t, \chi)| > e^{j (1/2-\delta_1) / D}$, so that, for any $\ell_{j,v} \geq 1$, we have $|Q_{j, 1}(t, \chi)|^{2\ell_{j, v}} e^{-2 j (1/2-\delta_1)\ell_{j, v} / D} \geq 1$. Therefore, by subdividing the range of integration $\mathcal{T}_2(\chi)$ into sets according to this $j$, we get
\[
\begin{split}
I_2 &\ll D \log Q_2 \sum_{\substack{\lfloor D \log P_2 \rfloor \leq v \leq D \log Q_2 \\ \lfloor D \log P_1 \rfloor \leq j \leq D \log Q_1}} e^{2(1/2- \delta_2) v / D - 2(1/2-\delta_1) \ell_{j, v} j /D} \\
&\qquad \qquad \cdot \sum_{\chi \pmod{q_1}} \int_\mathcal{T} |Q_{j, 1}(t, \chi)^{\ell_{j, v}} R_{v, 2}(t, \chi)|^2 dt
\end{split}
\]
where we pick for each $j,v$ an integer $\ell_{j,v}$ such that $\ell_{j,v} = \log (e^{v/D}) / \log (e^{j / D}) + O(1) = v / j + O(1)$. At this point we save so much that we do not need to exploit the $q$ average, so we can apply~\cite[Lemma 13]{mr} using point-wise bounds $k^{2\ell_{j, v}}$ and $(\log X)^{\nopf k \log k}$ for the coefficients of $Q_{j,1}^{2\ell_{j, v}}$ and $R_{v, 2}$ getting
\[
\begin{split}
\int_{\mathcal{T}} |Q_{j, 1}(t, \chi)^{\ell_{j, v}} R_{v, 2}(t, \chi)|^2 dt &\ll X k^{2\ell_{j, v}} (\log X)^{2\nopf k \log k + 10} \left(\frac{1}{qQ} + \frac{1}{H} +2^{\ell_{j, v}} Q_1\right) (\ell_{j, v} +1)!^2 \\
&\ll X \exp((\log\log X)^{3/2+1/100}) 
\end{split}
\]
And so we conclude that
$$
I_2 \ll q D^3 (\log Q_2)^3 \exp(- \tfrac{1}{100} (\log\log X)^2 + 2(\log\log X)^{3/2+1/100}) 
\ll \exp( - \tfrac{1}{200} (\log\log X)^{2})
$$
which is of course sufficient. 

Finally it remains to bound 
$$
I_3 = D \log Q_3 \sum_{\lfloor D \log P_3 \rfloor \leq v \leq D \log Q_3} \sum_{\chi \pmod{q_1}} \int_{\mathcal{T}_3(\chi)} |Q_{v,3}(t, \chi) R_{v, 3}(t, \chi)|^2 dt.
$$
We can replace the integral over $\mathcal{T}_3(\chi)$ by a sum over a one-spaced subset $\mathcal{U}(\chi) \subset \mathcal{T}_3(\chi)$. For each $t \in \mathcal{U}(\chi) \subset \mathcal{T}_3(\chi)$ there exists a $\lfloor D \log P_2 \rfloor \leq v \leq D \log Q_2$ such that $|Q_{v,2}(t, \chi)| > e^{(1/2-\delta_2) v / D}$. This implies by~\cite[Lemma 8]{mr} that $|\mathcal{U}(\chi)| \ll X^{1/2 - 1/1000}$.

Since the Dirichlet polynomial $Q_{v,3}(t, \chi)$ has length at least $P_3$, the Vinogradov-Korobov zero-free region for $L(s, \chi)$ implies that, for every $t \in \mathcal{T}$, we have $|Q_{v,3}(t, \chi)| \ll e^{v/(2D)} (\log X)^{10}/Q^3$. Therefore
$$
I_3 \ll \frac{1}{Q^4} \sup_{\substack{\lfloor D \log P_3 \rfloor \leq v \leq D \log Q_3 \\ \chi \pmod{q_1}}} e^{v/D} \sum_{t\in \mathcal{U}(\chi)} |R_{v,3}(t, \chi)|^2.
$$
Applying the Montgomery-Halasz lemma (see \cite[Theorem 9.6]{ik}) we see that 
\[
\begin{split}
I_3 &\ll \frac{1}{Q^4} \sup_{\substack{\lfloor D \log P_3 \rfloor \leq v \leq D \log Q_3 \\ \chi \pmod{q_1}}} e^{v/D} (X/e^{v/D} + |\mathcal{U}(\chi)| X^{1/2}) \sum_{n \sim X/(q_0e^{v/D})} \frac{d_k(n)^2}{n} \ll \frac{X}{Q^3}.
\end{split}
\]
which is sufficient.

\section{Large value estimates}\label{lve-sec}

In this section we establish Proposition \ref{lve-prop}.  We focus on the more difficult part (i) of this proposition, and indicate at the end of the section the changes needed to handle the (technically simpler) part (ii).  Thus in the following discussion we set $f \colon \Z \to \C$ to be the function
$$ f(n) \coloneqq d_k(n) 1_{{\mathcal S}_{k,X}}(n).$$ Note that by the definition of ${\mathcal S}_{k, X}$, we have 
\begin{equation}
\label{eq:f(n)bound}
f(n) \leq k^{\nopf k \log \log X} = (\log X)^{\nopf k \log k}.
\end{equation}

We first dispose of an easy case when $J$ is very large.  From the $1/H$-separated nature of the $\alpha_j$, Plancherel's identity,~\eqref{eq:f(n)bound} and divisor bound~\eqref{divisor-bound} we have
\begin{align*}
\sum_{j=1}^J \int_{[\alpha_j-1/2H,\alpha_j+1/2H]} |S_f(\beta;X)|^2\ d\beta 
&\leq \int_{\T} |S_f(\beta;X)|^2\ d\beta \\
&= \sum_n |f(n)|^2 \\
&\leq (\log X)^{\nopf k \log k} \sum_n |f(n)| \\
&\ll_k X \log^{\nopf k \log k + k-1} X.
\end{align*}
This bound gives the desired estimate \eqref{fx} unless we are in the case
\begin{equation}\label{J-bound}
J \leq \log^{5 k \log k} X,
\end{equation}
which we will now assume henceforth.

We first need to majorize the function $f = d_k 1_{{\mathcal S}_{k,X}}$ by a more tractable divisor sum $\tilde d_k$, defined by
\begin{equation}\label{tdk-def}
 \tilde d_k(n) \coloneqq \sum_{\substack{m|n: m \leq M \\ \Omega(m) \leq \nopf k\log \log X}} d_{k-1}(m) 1_{(X,2X]}(n)
\end{equation}
where $M$ is the quantity
\begin{equation}\label{M-def}
M \coloneqq X^{\frac{2k}{\log\log X}}.
\end{equation}
Clearly we have the upper bound
\begin{equation}\label{tkdk}
\tilde d_k(n) \leq d_k(n) 1_{(X,2X]}.
\end{equation}
In the opposite direction, we can use the construction of ${\mathcal S}_{k,X}$ to obtain

\begin{lemma}[Majorant property]\label{ndic}  For any integer $n$, one has the upper bound
$$ f(n) \ll \tilde d_k(n) (\log \log X)^{O_k(1)}.$$
\end{lemma}

One could remove the $(\log \log X)^{O(1)}$ loss here by iterating the construction of our majorant on the larger primes as done for example in the Brun-Hooley sieve \cite{fordhooley} or by replacing the majorant $\tilde d_k$ by the closely related majorants of Matthiesen \cite{matt, matt2}. In any case we will not need to do so here as we will (barely) be able to tolerate losses of the form $\log^{o(1)} X$ in our arguments.

\begin{proof}
For $n \in {\mathcal S}_{k, X}$, let us write $n = n_1 n_2$ where all prime factors of $n_1$ are at most $X^{1/(\log \log X)^2}$ and all prime factors of $n_2$ are larger than $X^{1/(\log \log X)^2}$. By definition of ${\mathcal S}_{k, X}$, we have $n_1 \leq X^{\nopf k \log \log X/(\log \log X)^2} \leq M$ and
\[
d_{k}(n) = d_{k}(n_1) d_{k}(n_2) \leq d_{k}(n_1) k^{10k \log \log \log X} \leq d_{k}(n_1) (\log \log X)^{O_k(1)}.
\]
Hence
\begin{align*}
f(n) &\leq 1_{S_{k, X}}(n)(\log \log X)^{O_k(1)} d_{k}(n_1) \\
&\leq 1_{(X, 2X]}(n) (\log \log X)^{O_k(1)}  \sum_{\substack{m \mid n_1 \\ \Omega(m)\\ \leq \nopf k\log \log x}} d_{k-1}(m) \\
&\leq 1_{(X, 2X]}(n) (\log \log X)^{O_k(1)} \sum_{\substack{m \mid n \\ m \leq M \\ \Omega(m) \leq \nopf k\log \log x}} d_{k-1}(m) \\
&= \tilde{d}_k(n) (\log \log X)^{O_k(1)}.
\end{align*}
\end{proof}

Now, we study the behavior of $\tilde d_k$ in short intervals and arithmetic progressions. First we have the following Brun-Titchmarsh type bounds on average.  

\begin{proposition}\label{excep}  Let $X \geq H_1 \geq (\log X)^{(1+2\varepsilon')k \log k}$. Let further
\[
q \leq \min\{H_1/(\log X)^{(1+2\varepsilon') k \log k}, (\log X)^{10000k^{10}}\}
\]
and $(a, q) = 1$. Then there exists an exceptional set ${\mathcal E} \subset [X/2,2X]$ (depending on $a, q, H_1, X, k$) of Lebesgue measure
\[
m({\mathcal E}) \ll_k X \log^{-10000 k^{10}} X + X \frac{(\log X)^{3(1+\varepsilon')k \log k+7}}{(H_1/q)^2}
\]
such that for all $x \in [X/2,2X] \backslash {\mathcal E}$, one has
\begin{equation}\label{jup2}
\sum_{\substack{x \leq n \leq x + H_1 \\ n = a\ (q)}} \tilde d_k(n) \ll_k \frac{H_1}{q} \log^{k-1+o(1)} X.
\end{equation}
\end{proposition}

\begin{proof} 
Let $x \in [X/2,3X]$ and let $l \geq 1$ be a fixed integer. Observe that for any finite sequence $a_1 \geq \dots \geq a_r$ of non-negative reals, one has the estimate
\begin{align*}
(\sum_{i=1}^r a_i)^l &\ll_l  \sum_{1 \leq i_1 \leq \dots \leq i_l \leq r} a_{i_1} \dots a_{i_l} \\
&\leq \sum_{1 \leq i_1 < \dots < i_l \leq r} a_{i_1} \dots a_{i_l} + (l-1) \left(\sup_{1 \leq i \leq r} a_i\right) \sum_{1 \leq i_1 \leq \dots \leq i_{l-1} \leq r} a_{i_1} \dots a_{i_{l-1}} \\
&\ll_l \sum_{1 \leq i_1 < \dots < i_l \leq r} a_{i_1} \dots a_{i_l} + \left(\sup_{1 \leq i \leq r} a_i\right) (\sum_{i=1}^r a_i)^{l-1}
\end{align*}
and hence
\begin{equation}
\label{eq:SupCorBound} 
(\sum_{i=1}^r a_i)^l \ll_l \sum_{1 \leq i_1 < \dots < i_l \leq r} a_{i_1} \dots a_{i_l}  + \left(\sup_{1 \leq i \leq r} a_i\right)^l.
\end{equation}
In particular
\begin{equation}
\label{eq:dklmoment}
\left(\sum_{\substack{x \leq n \leq x+H_1:\\ n = a\ (q)}} \tilde d_k(n)\right)^l \ll_l \left(\sup_{\substack{x \leq n \leq x+H_1 \\ n = a\ (q)}} \tilde d_k(n)\right)^l + \sum_{\substack{x \leq n_1 < \dots < n_l \leq x+H_1 \\ n_1,\dots,n_l = a\ (q)}} d_k(n_1) \dots d_k(n_l).
\end{equation}

Let us first dispose of exceptionally large values of the sup-term. Write
\[
\mathcal{N} = \left\{n \leq 4X \colon \tilde d_k(n) \geq \frac{H_1}{q} (\log X)^{k-1} \right\}.
\]
Now 
\[
\begin{split}
\# \mathcal{N} \cdot \left(\frac{H_1}{q}\right)^3 (\log X)^{3(k-1)} &\leq \sum_{n \leq 4X} \tilde d_k(n)^3 \\
&= \sum_{\substack{m_1, m_2, m_3 \leq M \\ \Omega(m_i) \leq (1+\varepsilon')k \log \log X}} d_{k-1}(m_1) d_{k-1}(m_2) d_{k-1}(m_3) \sum_{\substack{n \leq 4X \\ [m_1, m_2, m_3] \mid n}} 1 \\
&\ll X (\log X)^{3 \nopf k \log k} \sum_{\substack{m_1, m_2, m_3 \leq M}} \frac{1}{[m_1, m_2, m_3]} \\
&\ll X (\log X)^{3 \nopf k \log k + 7}.
\end{split}
\]
Hence 
\[
\#\mathcal{N} \ll \frac{X (\log X)^{3 \nopf k \log k + 7}}{\left(\frac{H_1}{q}\right)^3 (\log X)^{3(k-1)}}.
\]
Writing
\[
\mathcal{E}_1 = \left\{X/2 \leq x \leq 4X \colon \sup_{\substack{n \in [x, x+H_1] \\ n \equiv a \pmod{q}}} \tilde d_k(n) \geq \frac{H_1}{q} (\log X)^{k-1}\right\},
\]
we have $m(\mathcal{E}_1) \ll \#\mathcal{N} \cdot H_1/q$ and by the bound for $\# \mathcal{N}$ the set $\mathcal{E}_1$ can be included in the exceptional set $\mathcal{E}$. Hence it remains to show that \eqref{jup2} holds for all $x \in [X/2, 3X] \setminus \mathcal{E}_1$ with acceptably many exceptions.

Integrating \eqref{eq:dklmoment} over $[X/2,3X] \setminus \mathcal{E}_1$, we conclude that
\begin{align*}
&\int_{[X/2,3X] \setminus \mathcal{E}_1} \left(\sum_{\substack{x \leq n \leq x+H_1\\ n = a\ (q)}} \tilde d_k(n)\right)^l dx \\
&\ll_l \int_{[X/2,3X] \setminus \mathcal{E}_1} \left(\sup_{\substack{x \leq n \leq x+H_1 \\ n = a\ (q)}} \tilde d_k(n)\right)^l dx + \int_{X/2}^{3X} \sum_{\substack{y \leq n_1 < \dots < n_l \leq y+H_1 \\ n_1,\dots,n_l = a\ (q)}} d_k(n_1) \dots d_k(n_l) dy \\
&\ll_l X \left(\frac{H_1}{q}\right)^l (\log X)^{l(k-1)} + \sum_{0 = h_1 < h_2 < \dots < h_l \leq H_1/q} H_1 \sum_{n \leq 4X/q} d_k(qn+qh_1+a) \dots d_k(qn+qh_l+a).
\end{align*}

Using Lemma \ref{henriot}(ii) and the distinct nature of the $h_1,\dots,h_l$, we have
$$ \sum_{n \leq 4X/q} \tilde d_k(qn+qh_1+a) \dots \tilde d_k(qn+qh_l+a) \ll_{k,l} \frac{X}{q} \log^{l(k-1)+o(1)} X$$
On summing, we conclude that
$$ \sum_{0=h_1 < h_2 < \dots < h_l \leq H_1/q} H_1 \sum_{n \leq 4X/q} \tilde d_k(qn+qh_1+a) \dots \tilde d_k(qn+qh_l + a) \ll_{k,l} X \left(\frac{H_1}{q}\right)^l \log^{l(k-1)+o(1)} X.$$

Combining the bounds, we conclude that
$$
\int_{[X/2,3X] \setminus \mathcal{E}_1} \left(\sum_{\substack{x \leq n \leq x+H_1\\ n = a\ (q)}} \tilde d_k(n)\right)^l dx \ll_{k,l} X \left(\frac{H_1}{q}\right)^l \log^{l(k-1)+o(1)} X.$$
By Chebyshev's inequality, we conclude that
$$ 
\sum_{\substack{x \leq n \leq x+H_1\\ n = a\ (q)}} \tilde d_k(n) \leq \frac{H_1}{q} \log^{k-1 + o(1) + \frac{10000 k^{10}}{l}} X $$
for all $x \in [X/2,3X] \setminus \mathcal{E}_1$ outside of a set of measure at most $X \log^{-10000 k^{10}} X$, for any fixed large enough $l$.  Sending $l$ sufficiently slowly to infinity, the claim follows.
\end{proof}

For any $\alpha \in \T$ and $Q \geq 1$, let $q_{\alpha,Q}$ be the least positive integer $q$ such that
$$ \| \alpha - \frac{a}{q}\|_{\T} \leq \frac{1}{qQ}$$
for some $0 \leq a < q$ coprime to $q$, where $\|x\|_{\T}$ denotes the distance of $x$ to the integers; such a $q$ exists and is bounded by $Q$ thanks to the Dirichlet approximation theorem.  We now have the following variant of the above proposition:

\begin{proposition}\label{doc0} Let $X \geq H_1 \geq (\log X)^{10 k \log k}$, and let 
\begin{equation}
\label{eq:Q1bound}
Q_1 \leq \min\{H_1^{2/5}/(\log X)^{(1+2\varepsilon)k \log k}, (\log X)^{10000k^{10}} \}.
\end{equation}
Let $\alpha \in \T$ with $|\alpha - a/q| < \frac{1}{qQ_1}$ where $q \leq Q_1$ and $(a, q) = 1$.  There exists an exceptional set ${\mathcal E} \subset [X/2,2X]$ (possibly depending on $\alpha$) of Lebesgue measure
$$ m({\mathcal E}) \ll_k X(\log X)^{3k \log k}/Q_1^{1/4}$$ such that for all $x \in [X/2,2X] \backslash {\mathcal E}$, one has
\begin{equation}
\label{eq:doc0claim}
\sum_{x \leq n \leq x+H_1} \tilde d_k(n) e(\alpha n) \ll_k \left(\frac{1}{q} + \frac{1}{Q_1^{1/8}}\right) H_1 \log^{k-1+o(1)} X.
\end{equation}
\end{proposition}

\begin{proof} Write $\alpha = \frac{a}{q} + \frac{\theta}{qQ_1}$ for some $\theta \in [-1,1]$. 
Splitting the $n$ summation in~\eqref{eq:doc0claim} according to the value of $q' = (n, q)$, we can bound the left hand side of~\eqref{eq:doc0claim} by 
\begin{equation}
\label{eq:addtomultchars}
\begin{split}
&\sum_{q = q' q_2} \Bigl|\sum_{\substack{\frac{x}{q'} \leq n \leq \frac{x}{q'} + \frac{H_1}{q'}}} \tilde d_{k}(q'n) e\left(\frac{an}{q_2}\right) 1_{(n,q_2)=1} e\left(\frac{\theta n}{q_2 Q_1}\right)\Bigr| \\
&\leq \sum_{q = q' q_2} \frac{1}{\varphi(q_2)} \sum_{\chi \pmod{q_2}} \Bigl|\tau(\overline{\chi}) \sum_{\substack{\frac{x}{q'} \leq n \leq \frac{x}{q'} + \frac{H_1}{q'}}} \tilde d_{k}(q'n) \chi(n) e\left(\frac{\theta n}{q_2 Q_1}\right)\Bigr|,
\end{split}
\end{equation}
where we have decomposed $e(\frac{an}{q_2}) 1_{(n,q_2)=1}$ into Dirichlet characters
$$ 
e\left(\frac{an}{q_2}\right) 1_{(n,q_2)=1} = \frac{1}{\phi(q_2)} \sum_{\chi\ (q_2)} \chi(a) \chi(n) \tau(\overline{\chi})
$$
where $\tau(\overline{\chi})$ is the Gauss sum
$$ \tau(\overline{\chi}) \coloneqq \sum_{l=1}^{q_2} e\left(\frac{l}{q_2}\right) \overline{\chi(l)}.$$
We first consider the contribution when $\chi$ is a principal character.  In this case $\tau(\overline{\chi})$ is a Ramanujan sum and so $|\tau(\overline{\chi})| \leq 1$ (see e.g.~\cite[formula (3.4)]{ik}). Hence the contribution of the principal character to \eqref{eq:addtomultchars} is at most
\[
\begin{split}
&d(q) \sup_{q=q' q_2} \frac{1}{\varphi(q_2)} \left|\sum_{\substack{\frac{x}{q'} \leq n \leq \frac{x}{q'} + \frac{H_1}{q'}}} \tilde d_k(q'n) \right| \leq (\log X)^{o(1)} \sup_{q = q' q_2} \frac{1}{q_2} \left|\sum_{\substack{x \leq n \leq x + H_1 \\ n \equiv 0 \pmod{q'}}} \tilde d_{k}(n) \right|.
\end{split}
\]
By Proposition~\ref{excep} the total contribution of the principal characters is acceptable due to our conditions on sizes of $H_1$ and $Q_1$.

For the non-principal characters, we bound $|\tau(\overline{\chi})| \leq \sqrt{q_2}$ and use~\eqref{tdk-def} to bound their contribution to~(\ref{eq:addtomultchars}) by 
\[
\sum_{q = q' q_2} \frac{\sqrt{q_2}}{\varphi(q_2)} \sum_{\substack{\chi \pmod{q_2} \\ \chi \neq \chi_0}} \left| \sum_{\substack{x \leq mr \leq x + H_1 \\ q' \mid mr, m \leq M \\ \Omega(m) \leq \nopf k\log \log x }} d_{k-1}(m) 1_{(X,2X]}(mr) \chi(mr/q') e\left(\frac{\theta mr}{q'q_2 Q_1}\right)\right|.
\]
We split the summations over $m$ and $r$ according to $q_0 \mid q'$ such that $q_0 \mid m, q_1 = \frac{q'}{q_0} \mid r$ and $(m/q_0, q_2) = (r/q_1, q_2) = 1$ (choice of $q_0$ might not be unique, but it always exists). Then the above expression is at most
\[
\begin{split}
&\sum_{q = q_0 q_1 q_2} \frac{\sqrt{q_2}}{\varphi(q_2)} \sum_{\substack{\chi \pmod{q_2} \\ \chi \neq \chi_0}} \left|\sum_{\substack{\frac{x}{q_0 q_1} \leq mr \leq \frac{x}{q_0 q_1} + \frac{H_1}{q_0 q_1} \\ mq_0 \leq M \\ \Omega(mq_0) \leq \nopf k\log \log x}} d_{k-1}(mq_0) \chi(mr) e\left(\frac{\theta mr}{q_2 Q_1}\right)\right| dx.
\end{split}
\]

By Cauchy-Schwarz and the triangle inequality, the claim follows once we have shown for all $q=q_0q_1q_2$ and all non-principal characters $\chi$ of period $q_2$ that
$$
\int_{X}^{2X} \left|\sum_{\substack{\frac{x}{q_0 q_1} \leq mr \leq \frac{x}{q_0 q_1} + \frac{H_1}{q_0 q_1} \\ mq_0 \leq M \\ \Omega(mq_0) \leq \nopf k\log \log x}} d_{k-1}(q_0 m) \chi(mr) e\left(\frac{\theta mr}{q_2 Q_1}\right)\right|^2 dx \ll_k \frac{H_1^2 X \log^{3k\log k+2k-2+o(1)} X}{q_2 Q_1^{1/2}}.
$$

Squaring out and writing $h = m_2r_2 - m_1 r_1$, the left hand side becomes
\[
\begin{split}
&\sum_h e\left(\frac{\theta h}{q_2 Q_1}\right) \sum_{\substack{m_i \leq M/q_0 \\ \Omega(m_i q_0) \leq \nopf k\log \log x}} \sum_{\substack{r_1, r_2 \\ m_2 r_2 - m_1 r_1 = h}} d_{k-1}(q_0 m_1) d_{k-1}(q_0m_2) \overline{\chi}(m_1 r_1) \chi(m_1 r_1+h)  \\
& \qquad \cdot m( (X,2X] \cap [q_0 q_1 m_1 r_1 - H_1, q_0 q_1 m_1 r_1] \cap [q_0 q_1 m_2 r_2 - H_1, q_0 q_1 m_2 r_2]).
\end{split}
\]
The quantity on the second line vanishes unless $|h| \leq \frac{H_1}{q_0 q_1}$ and $m_1 r_1 \in [\frac{X}{q_0q_1}, \frac{2X}{q_0q_1} + \frac{H_1}{q_0q_1}]$. If these two conditions hold, the quanity on the second line is equal to $H_1 - q_0 q_1 |h|$.  The contribution of those boundary $m_1r_1$ with $m_1 r_1 = \frac{X}{q_0q_1} + O(\frac{H_1}{q_0q_1})$ or $m_1 r_1 = \frac{2X}{q_0q_1} + O(\frac{H_1}{q_0q_1})$ is easily seen to be acceptable, so it suffices to show that
$$ \sum_{|h| \leq \frac{H_1}{q_0 q_1}} (H_1-q_0 q_1 |h|) e\left(\frac{\theta h}{q_2 Q_1}\right) A(h)
\ll_k \frac{H_1^2 X \log^{3k\log k+2k-2+o(1)} X}{q_2 Q_1^{1/2}}, $$
where $A(h)$ is the quantity
$$ A(h) \coloneqq \sum_{\substack{m_1,m_2 \leq M/q_0, (m_1 m_2, q_2) = 1 \\ \Omega(m_1q_0), \Omega(m_2q_0) \leq \nopf k \log \log x}} d_{k-1}(m_1q_0) d_{k-1}(m_2q_0) A_{m_1,m_2}(h)$$
where
$$ A_{m_1,m_2}(h) \coloneqq \sum_{\substack{n \in (\frac{X}{q_0 q_1},\frac{2X}{q_0 q_1}] \\ m_1 | n; m_2 | n+h}} \overline{\chi}(n) \chi(n+h).$$

Now
\[
A_{m_1,m_2}(h) = \sum_{b \pmod{q_2}} \overline{\chi}(b) \chi(b+h) \sum_{\substack{n \sim X/(q_0q_1) \\ n \equiv b \pmod{q_2} \\  m_1 \mid n \\ m_2 \mid n +h}} 1.
\]
The system of equations on $n$ is soluble only when $(m_1, m_2) \mid h$ in which case it has an unique solution $\pmod{q_2[m_1, m_2]}$. Hence, when these conditions hold
\[
A_{m_1,m_2}(h) = \frac{X}{q_0 q_1 q_2 [m_1, m_2]} \sum_{b \pmod {q_2}} \overline{\chi}(b) \chi(b+h) + O(q_2).
\]
The contribution of the $O(q_2)$ error is bounded by
$$ \sum_{|h| \ll H_1} H_1 \sum_{m_1,m_2 \leq M/q_0} d_{k-1}(m_1q_0) d_{k-1}(m_2q_0) Q_1 $$
which is easily seen to be acceptable from \eqref{M-def}, \eqref{eq:Q1bound}. Hence we end up with
\[
\begin{split}
&\frac{X}{q_0 q_1 q_2} \sum_{\substack{m_1 \leq M/q_0 \\ \Omega(m_1 q_0) \leq \nopf k \log \log x \\ (m_1, q_2) = 1}} \sum_{\substack{m_2 \leq M/q_0 \\ \Omega(m_2 q_0) \leq \nopf k \log \log x \\ (m_2, q_2) = 1}} \frac{d_{k-1}(q_0 m_1) d_{k-1}(q_0 m_2)}{[m_1, m_2]} \\ &\cdot \sum_{|h| \leq \frac{H_1}{q_0 q_1 (m_1, m_2)}} (H_1-q_0 q_1(m_1, m_2) |h|) e\left(\frac{\theta h (m_1, m_2)}{q_2 Q_1}\right) \sum_{b \pmod{q_2}} \overline{\chi}(b) \chi(b+h (m_1, m_2)).
\end{split}
\]
Since $(m_1m_2, q_2) = 1$, here 
\[
\sum_{b \pmod{q_2}} \overline{\chi}(b) \chi(b+h (m_1, m_2)) = \sum_{b \pmod{q_2}} \overline{\chi}(b) \chi(b+h)
\]
Thus, by the triangle inequality, it suffices to show that, for any $(b, q_2) = 1$,
\begin{equation}\label{dad}
\begin{split}
&\frac{X}{q_0 q_1} \sum_{\substack{m_1 \leq M/q_0 \\ \Omega(m_1 q_0) \leq \nopf k \log \log x \\ (m_1, q_2) = 1}} \sum_{\substack{m_2 \leq M/q_0 \\ \Omega(m_2 q_0) \leq \nopf k \log \log x \\ (m_2, q_2) = 1}} \frac{d_{k-1}(m_1 q_0) d_{k-1}(m_2 q_0)}{[m_1, m_2]} |B(m_1,m_2)| \\
&\quad \ll_k \frac{H_1^2 X \log^{3k\log k +2k-2+o(1)} X}{q_2 Q_1^{1/2}},
\end{split}.
\end{equation}
where
$$ B(m_1,m_2) \coloneqq \sum_{|h| \leq \frac{H_1}{q_0q_1(m_1, m_2)}} (H_1-q_0 q_1 (m_1, m_2)|h|) e\left(\frac{\theta h (m_1, m_2)}{q_2 Q_1}\right) \chi(b+h).$$
Partial summation and Polya-Vinogradov give
\[
\begin{split} 
B(m_1,m_2) &\ll \left(H_1 + \frac{\theta H_1^2}{q_0 q_1 q_2 Q_1} \right) \max_{x \leq H_1} \left|\sum_{h \leq x}\chi(b+h)\right| \\
&\ll H_1 q_2^{1/2}\log q_2 + \frac{H_1^2 \log q_2}{q_0 q_1 q_2^{1/2} Q_1}.
\end{split}
\]
Hence the left hand side of~\eqref{dad} is at most 
\begin{align*}
& \left( \frac{H_1 X q_2^{1/2} \log q_2}{q_0 q_1} + \frac{H_1^2 X}{(q_0 q_1)^2 q_2^{1/2} Q_1} \log q_2\right) \\ \times & \sum_{\substack{m_1 \leq M/q_0 \\ \Omega(m_1 q_0) \leq \nopf k \log \log x \\ (m_1, q_2) = 1}} \sum_{\substack{m_2 \leq M/q_0 \\ \Omega(m_2 q_0) \leq \nopf k \log \log x \\ (m_2, q_2) = 1}} \frac{d_{k-1}(q_0 m_1) d_{k-1}(q_0 m_2)}{[m_1, m_2]}. 
\end{align*}
Writing $d = (m_1, m_2)$, the sum over $m_1$ and $m_2$ is at most
\[
\begin{split}
&\sum_{\substack{d \leq M \\ \Omega(d) \leq \nopf k \log \log x}} \frac{d_{k-1}(d)^2}{d} \sum_{\substack{m_1 \leq M/d \\ \Omega(m_1) \leq \nopf k \log \log x}} \sum_{\substack{m_2 \leq M/d \\ \Omega(m_2) \leq \nopf k \log \log x}} \frac{d_{k-1}(m_1) d_{k-1}(m_2)}{m_1 m_2} \\
&\ll (\log X)^{\nopf k\log k + 3k-3},  
\end{split}
\]
and the claim follows.
\end{proof}

We will actually use the following consequence of Proposition~\ref{doc0}.
\begin{proposition}\label{doc} Let $X \geq H_1 \geq (\log X)^{20 k \log k}$, and let 
\[
Q_1 \leq \min\{H_1^{1/3}/(\log X)^{k \log k}, (\log X)^{10 000k^{10}} \}.
\]
Let $\alpha \in \T$ with $|\alpha - a/q| < \frac{1}{qQ_1}$ where $q \leq Q_1$ and $(a, q) = 1$.  There exists an exceptional set ${\mathcal E} \subset [X/2,2X]$ (possibly depending on $\alpha$) of Lebesgue measure
$$ m({\mathcal E}) \ll_k X (\log X)^{3k\log k}/Q_1^{1/8}$$ such that for all $x \in [X/2,2X] \backslash {\mathcal E}$, one has
$$ \sup_{0 \leq H' \leq H_1} \left|\sum_{x \leq n \leq x+H'} \tilde d_k(n) e(\alpha n)\right| \ll_k \left(\frac{1}{q} + \frac{1}{Q_1^{1/8}}\right) H_1 \log^{k-1+o(1)} X.$$
\end{proposition}

\begin{proof} We round $H'$ to the nearst multiple of $H_1/Q_1^{1/8}$ and apply Proposition~\ref{doc0} with $H_1$ there $H_1/Q_1^{1/8}$ and use the union bound.
\end{proof}

We return to the proof of Proposition \ref{lve-prop}(i).  It will suffice to establish the bound
$$
\sum_{j=1}^J \left(\int_{[\alpha_j-1/2H,\alpha_j+1/2H]} |S_f(\beta;X)|^2\ d\beta\right)^{1/2} \ll_k J^{5/8} X^{1/2} \log^{k-1+o(1)} X$$
since (after discarding some of the $\alpha_j$ as necessary) this implies that for any $1 \leq J' \leq J$, the $(J')^{\operatorname{th}}$ largest value of $(\sum_{[\alpha_j-1/2H,\alpha_j+1/2H]} |S_f(\beta;X)|^2\ d\beta)^{1/2}$ is $O_k( (J')^{-3/8} X^{1/2} \log^{k-1+o(1)} X)$, and the claim then follows by square-summing in $J'$.

By Gallagher's lemma (Lemma~\ref{le:Gallagher}), we have
$$ \left(\int_{[\alpha_j-1/2H,\alpha_j+1/2H]} |S_f(\beta;X)|^2\ d\beta\right)^{1/2} \ll \frac{1}{H} \left(\int_\R |\sum_{x \leq n < x+\frac{H}{10}} f(n) e(\alpha_j n)|^2\ dx\right)^{1/2},$$
and so it suffices to show that
$$ \sum_{j=1}^J \left(\int_\R |\sum_{x \leq n < x+\frac{H}{10}} f(n) e(\alpha_j n)|^2\ dx\right)^{1/2} \ll_k J^{5/8} X^{1/2} H \log^{k-1+o(1)} X.$$
We may restrict the $x$ variable to the interval $[X/2,2X]$, since the integral vanishes otherwise. Moreover, we can assume that for each $j$ we have
\begin{equation}
\label{eq:intjlow}
\int_{X/2}^{2X} \left|\sum_{x \leq n < x+\frac{H}{10}} f(n) e(\alpha_j n)\right|^2\ dx \geq \frac{1}{J} X H^2 \log^{2k-2} X.
\end{equation}
since those $j$ for which this does not hold make an acceptable contribution and can be discarded.
We will need the quantity
\begin{equation}\label{Q-def}
Q \coloneqq \log^{3000k \log k} X;
\end{equation}
this is a parameter that is much larger than $J$ (recall~(\ref{J-bound})), but much smaller than $H$.  Applying Proposition~\ref{doc} to each $\alpha_i - \alpha_j$ with $(H_1, Q_1) = (H, Q)$ and to each $a/q$, $q \leq J^2$ and $(a, q) =1$ with $(H_1, Q_1) =  (Q^{3/4}, Q^{1/4})$, and taking unions, we can find an exceptional set ${\mathcal E} \subset [X/2,2X]$ of Lebesgue measure
\begin{equation}
\label{eq:measE''}
m({\mathcal E}) \ll_k J^{4} (\log X)^{3k\log k} X/Q^{1/32} \ll \frac{X}{(\log X)^{70 k\log k}}
\end{equation}
(where we have recalled~(\ref{J-bound})) such that, for all $x \in [X/2,2X] \backslash {\mathcal E}$, one has
\begin{equation}\label{flag-1}
 \sup_{0 \leq H' \leq H} \left|\sum_{x \leq n \leq x+H'} \tilde d_k(n) e((\alpha_j - \alpha_{j'})n)\right| \ll_k \left(\frac{H}{q_{\alpha_j-\alpha_{j'},Q}} + \frac{H}{Q^{1/8}}\right) \log^{k-1+o(1)} X
\end{equation}
for all $1 \leq j,j' \leq J$, and also
\begin{equation}\label{flag-2} 
\sup_{0 \leq H' \leq Q^{3/4}} \left|\sum_{x \leq n \leq x+H'} \tilde d_k(n) e(an/q)\right| \ll_k \frac{Q^{3/4}}{q} \log^{k-1+o(1)} X.
\end{equation}
whenever $q \leq J^{2}$ and $a$ is coprime to $q$.

We first consider the contribution of the exceptional set ${\mathcal E}$.  We may use~\eqref{eq:SupCorBound}, \eqref{eq:f(n)bound} and Lemma~\ref{henriot}(iii) to bound, for any $\theta \in \mathbb{T}$,
\begin{equation}
\begin{split}
\label{eq:4thmom}
&\int_{[X/2,2X]} \left|\sum_{x \leq n < x+\frac{H}{10}} f(n) e(\theta n)\right|^4\ dx \\
&\ll H\sum_{n \in [X/2,2X]} f(n)^4 + H \sum_{0 < h_1 < h_2 < h_3 \leq H/10} \sum_{n \in [X, 2X]} d_k(n) d_k(n+h_1) d_k(n+h_2) d_k(n+h_3)  \\
&\ll HX (\log X)^{4 \nopf k \log k} + H^4 X (\log X)^{4(k-1)} \\
&\ll H^4 X \log^{4(k-1)} X,
\end{split}
\end{equation}
and hence, by the Cauchy-Schwarz inequality and~(\ref{eq:measE''}),
\begin{equation}
\label{eq:E''contr}
\sum_{j=1}^J \left(\int_{{\mathcal E}} \left|\sum_{x \leq n < x+\frac{H}{10}} f(n) e(\alpha_j n)\right|^2\ dx \right)^{1/2} \ll J \left(\frac{1}{(\log X)^{35 k \log k}} H^2 X \log^{2(k-1)} X\right)^{1/2},
\end{equation}
which when combined with \eqref{J-bound} ensures that the contribution of the exceptional set is acceptable.  Thus it suffices to show that
$$ \sum_{j=1}^J \left(\int_{[X/2,2X] \backslash {\mathcal E}} \left|\sum_{x \leq n < x+\frac{H}{10}} f(n) e(\alpha_j n)\right|^2\ dx\right)^{1/2} \ll_k J^{5/8} X^{1/2} H \log^{k-1+o(1)} X.$$
By duality, it thus suffices to show that
\begin{equation}
\label{eq:DualClaim} 
\sum_{j=1}^J \int_\R \sum_{x \leq n < x+\frac{H}{10}} f(n) e(\alpha_j n) g_j(x) \ dx \ll_k J^{5/8} X^{1/2} H \log^{k-1+o(1)} X,
\end{equation}
where 
\[
g_j(x) \coloneqq 1_{x \in [X/2, 2X] \backslash {\mathcal E}} \frac{\overline{\sum_{x \leq n \leq x+H/10} f(n) e(\alpha_j n)}}{\left(\int_{[X/2,2X] \backslash {\mathcal E}} |\sum_{y \leq n < y+\frac{H}{10}} f(n) e(\alpha_j n)|^2\ dy\right)^{1/2}}.
\]
Here $g_j: \R \to \C$ are measurable functions supported on $[X/2,2X] \backslash {\mathcal E}$ with
\begin{equation}\label{gj2}
 \int_\R |g_j(x)|^2\ dx = 1
\end{equation}
for $1 \leq j \leq J$. Also by~\eqref{eq:intjlow},~(\ref{flag-1}) with $j=j'$, and~\eqref{eq:E''contr} we can assume that 
\begin{equation}
\label{eq:gjup}
|g_j(x)| \leq J^{1/2} (\log X)^{o(1)}/X^{1/2}
\end{equation} 
for all $x$ and $j$.

We rewrite the left-hand side of~\eqref{eq:DualClaim} as
$$ \sum_n f(n) \sum_{j=1}^J \int_{n-H/10}^n e(\alpha_j n) g_j(x)\ dx$$
and use Lemma \ref{ndic} to bound this by
$$ \sum_n \tilde d_k(n) \left|\sum_{j=1}^J \int_{n-H/10}^n e(\alpha_j n) g_j(x)\ dx\right| \log^{o(1)} X.$$
From \eqref{divisor-bound} and \eqref{tkdk} we have
$$ \sum_n \tilde d_k(n) \ll_k X \log^{k-1} X$$
and so by the Cauchy-Schwarz inequality, it will suffice to show that
\begin{equation}\label{temp}
 \sum_n \tilde d_k(n) \left|\sum_{j=1}^J \int_{n-H/10}^n e(\alpha_j n) g_j(x)\ dx\right|^2
\ll_k J^{5/4} H^2 \log^{k-1+o(1)} X.
\end{equation}
The left-hand side may be rearranged as
\begin{equation}\label{temp2}
 \sum_{1 \leq j,j' \leq J} \int_\R \int_\R g_j(x) \overline{g_{j'}(x')}
\sum_{n \in [x,x+H/10] \cap [x',x'+H/10]} \tilde d_k(n) e((\alpha_j-\alpha_{j'})n)\ dx dx'.
\end{equation}
One could attempt to control this sum purely using \eqref{flag-1}, but this turns out to be insufficient due to the fact that the $\alpha_j$ are only $1/H$-separated instead of $1/Q$-separated, so we could for instance have $q_{\alpha_{j}-\alpha'_j, Q} = 1$ for every pair $(j, j')$.  To rectify this, we use the greedy algorithm to find a $1/Q$-separated sequence $\beta_1,\dots,\beta_{J'}$ of elements of $\T$ for some $1 \leq J' \leq J$, such that each $\alpha_j$ is within $1/Q$ of at least one of the $\beta_i$, $i=1,\dots,J'$.  Thus we can find a partition
$$ \{1,\dots,J\} = \bigcup_{i=1}^{J'} {\mathcal J}_i$$
where for each $1 \leq i \leq J'$ and $j \in {\mathcal J}_i$, one has 
\begin{equation}\label{alphabet}
\| \alpha_j - \beta_i \|_\T \leq 1/Q.
\end{equation}

By the triangle inequality, it suffices to show that
\begin{equation}
\label{eq:Aii'claim}
\sum_{1 \leq i,i' \leq J'} |A_{i,i'}| \ll J^{5/4} H^2 (\log X)^{k-1+o(1)},
\end{equation}
where, for each $1 \leq i,i' \leq J'$, $A_{i,i'}$ denotes the quantity
$$
A_{i,i'} \coloneqq \sum_{j \in {\mathcal J}_i} \sum_{j' \in {\mathcal J}_{i'}} \int_\R \int_\R g_j(x) \overline{g_{j'}(x')}
\sum_{n \in [x,x+H/10] \cap [x',x'+H/10]} \tilde d_k(n) e((\alpha_j-\alpha_{j'})n)\ dx dx'.
$$
Let us call a pair $(i,i')$ \emph{good} if one has
\begin{equation}\label{qQ}
 q_{\alpha_j - \alpha_{j'}, Q} > J^2
\end{equation}
for all $j \in {\mathcal J}_i$ and $j' \in {\mathcal J}_{i'}$, and \emph{bad} otherwise.  

If $(i,i')$ is good, then we can use the triangle inequality to bound
$$ |A_{i,i'}| \leq \sum_{j \in {\mathcal J}_i} \sum_{j' \in {\mathcal J}_{i'}}  \int_\R \int_\R |g_j(x)| |g_{j'}(x')|
\left|\sum_{n \in [x,x+H/10] \cap [x',x'+H/10]} \tilde d_k(n) e((\alpha_j-\alpha_{j'})n)\right|\ dx dx'.$$
By construction, $|g_j(x)|$ is only non-zero when $x \not \in {\mathcal E}$, and the inner sum is only non-zero when $x'=x+O(H)$.  By \eqref{flag-1} and \eqref{qQ}, we have
$$ \sum_{n \in [x,x+H/10] \cap [x',x'+H/10]} \tilde d_k(n) e((\alpha_j-\alpha_{j'})n) \ll_k \frac{H}{J^2} \log^{k-1+o(1)} X,$$
while from Schur's test (or the elementary bound $|g_j(x)| |g_{j'}(x')| \leq |g_j(x)|^2 + |g_{j'}(x')|^2$) and \eqref{gj2} we have
$$ \int\int_{x' = x+O(H)} |g_j(x)| |g_{j'}(x')|\ dx dx' \ll H.$$
Moreover, since
$$
\sum_{(i,i') \text{ good}} \sum_{j \in \mathcal{J}_{i}} \sum_{j' \in \mathcal{J}_{i'}} 1 \ll J^2
$$
we get, 
$$ \sum_{(i, i') \text{ good}}|A_{i,i'}| \leq H J^2 \cdot \frac{H}{J^2} \log^{k-1+o(1)} X \leq H^2 \log^{k-1+o(1)} X,$$
so the contribution of the good $(i,i')$ to \eqref{eq:Aii'claim} is acceptable.

It remains to show that
\begin{equation}\label{ssss}
\sum_{(i,i') \text{ bad}} |A_{i,i'}| \ll J^{5/4} H^2 \log^{k-1+o(1)} X.
\end{equation}
Suppose that $(i,i')$ is bad, then there exist $j \in {\mathcal J}_i$ and $j' \in {\mathcal J}_{i'}$ such that
$$ q_{\alpha_j - \alpha_{j'}, Q} \leq J^2.$$
Thus there exists $q_{i,i'} \leq J^{2}$ and $a_{i,i'}$ coprime to $q_{i,i'}$ such that
$$ \left\| \alpha_j - \alpha_{j'} - \frac{a_{i,i'}}{q_{i,i'}} \right\|_\T \leq \frac{1}{q_{i,i'}Q} \leq \frac{1}{Q}$$
which in particular implies by~(\ref{alphabet}) that
\begin{equation}\label{ga}
 \left\|\beta_i - \beta_{i'} - \frac{a_{i,i'}}{q_{i,i'}} \right\|_\T \leq \frac{3}{Q}.
\end{equation}
We may rewrite $A_{i,i'}$ as
$$ A_{i,i'} = \int_\R \int_\R \sum_{n \in [x,x+H/10] \cap [x',x'+H/10]} \tilde d_k(n) e\left( \frac{a_{i,i'}n}{q_{i,i'}}\right) F_{i,i'}(n) G_{i,x}(n) \overline{G_{i',x'}(n)}\ dx dx'$$
where
\begin{align*}
F_{i,i'}(n) & \coloneqq e\left( (\beta_i - \beta_{i'} - \frac{a_{i,i'}}{q_{i,i'}}) n \right) \\
G_{i,x}(n) & \coloneqq  \sum_{j \in {\mathcal J}_i} g_j(x) e\left( (\alpha_j - \beta_i) n \right) \\
G_{i',x'}(n) & \coloneqq  \sum_{j \in {\mathcal J}_{i'}} g_j(x') e\left( (\alpha_j - \beta_{i'}) n \right).
\end{align*}
We can of course restrict the integral to the region $|x'-x| \leq H/10$, since the sum vanishes otherwise.
From Cauchy-Schwarz we have the crude upper bounds
\begin{align*}
|F_{i,i'}(n)| &\leq 1 \\
|G_{i,x}(n)| &\leq |\mathcal{J}_i|^{1/2} (\sum_{j=1}^J |g_j(x)|^2)^{1/2} \\
|G_{i',x'}(n)| &\leq |\mathcal{J}_{i'}|^{1/2} (\sum_{j=1}^J |g_j(x')|^2)^{1/2}
\end{align*}
and where of course $|\mathcal{J}_i| \leq J$ and $|\mathcal{J}_{i'}| \leq J$. 

Whenever $|n-n'| \leq Q^{3/4}$, \eqref{alphabet} and \eqref{ga} give the bounds
\begin{align*}
 |F_{i,i'}(n) - F_{i,i'}(n')| &\ll Q^{-1/4} \\
 |G_{i,x}(n) - G_{i,x}(n')| &\ll J^{1/2} Q^{-1/4} (\sum_{j=1}^J |g_j(x)|^2)^{1/2} \\
 |G_{i',x'}(n) - G_{i',x'}(n')| &\ll J^{1/2} Q^{-1/4} (\sum_{j=1}^J |g_j(x')|^2)^{1/2}.
\end{align*}
Hence
\begin{align*}
F_{i,i'}(n') G_{i,x}(n') \overline{G_{i',x'}(n')} &= F_{i,i'}(n) G_{i,x}(n) \overline{G_{i',x'}(n)} \\
&\quad + O\left( J Q^{-1/4} \left(\sum_{j=1}^J |g_j(x)|^2\right)^{1/2} \left(\sum_{j=1}^J |g_j(x')|^2\right)^{1/2} \right)
\end{align*}
whenever $|n-n'| \leq Q^{3/4}$.  Thus, for any interval $I_j \coloneqq [jQ^{3/4}, (j+1)Q^{3/4}] \subset [0,H/10]$ and any $n_{I_j} \in I_j$, we can write
\[
\begin{split}
&\sum_{n \in (x+I_j) \cap [x',x'+H/10]} \tilde d_k(n) e\left( \frac{a_{i,i'}n}{q_{i,i'}} \right) F_{i,i'}(n) G_{i,x}(n) \overline{G_{i',x'}(n)} \\
& = \sum_{n \in (x+I_j) \cap [x',x'+H/10]} \tilde d_k(n) e\left( \frac{a_{i,i'}n}{q_{i,i'}} \right) F_{i,i'}(n_{I_j}) G_{i,x}(n_{I_j}) \overline{G_{i',x'}(n_{I_j})}\\
&\quad + O\left( J Q^{-1/4} \left(\sum_{j=1}^J |g_j(x)|^2\right)^{1/2} \left(\sum_{j=1}^J |g_j(x')|^2\right)^{1/2} \sum_{n \in x+{I_j}} \tilde d_k(n) \right).
\end{split}
\]
We can apply \eqref{flag-2} and \eqref{eq:gjup} to bound this by
\[ 
\begin{split}
&\ll \frac{Q^{3/4}}{q_{i,i'}} \log^{k-1+o(1)} X |G_{i,x}(n_{I_j})| |G_{i',x'}(n_{I_j})| + 1_{x + j Q^{3/4} \in \mathcal{E}} \frac{J^3}{X} (\log X)^{o(1)} \sum_{n \in x+I_j} \tilde{d}_k(n) \\
&\qquad +  J Q^{1/2} \left(\sum_{j=1}^J |g_j(x)|^2\right)^{1/2} \left(\sum_{j=1}^J |g_j(x')|^2\right)^{1/2} \log^{k-1+o(1)} X,
\end{split}
\]
so on averaging over $n_{I_j}$ we obtain the bound
\[
\begin{split}
&\ll \frac{1}{q_{i,i'}} \log^{k-1+o(1)} X \sum_{n \in x+{I_j}} |G_{i,x}(n)| |G_{i',x'}(n)| + 1_{x + j Q^{3/4} \in \mathcal{E}} \frac{J^3}{X} (\log X)^{o(1)} \sum_{n \in x+I_j} \tilde{d}_k(n) \\ 
& \qquad +  J Q^{1/2} \left(\sum_{j=1}^J |g_j(x)|^2\right)^{1/2} \left(\sum_{j=1}^J |g_j(x')|^2\right)^{1/2} \log^{k-1+o(1)} X.
\end{split}
\]

The second term contributes to $A_{i, i'}$ at most
\[
\begin{split}
&H \frac{J^3}{X} (\log X)^{o(1)}\sum_{j=1}^{H/(10Q^{3/4})} \int_{X/2}^{2X} 1_{x + j Q^{3/4} \in \mathcal{E}} \sum_{n \in x+I_j} \tilde{d}_k(n) dx \\
&\ll H^2 \frac{J^3}{X} (\log X)^{o(1)} \int_{\mathcal{E}} \frac{1}{Q^{3/4}} \sum_{x \leq n < x + Q^{3/4}} \tilde{d}_k(n) dx \\
&\ll H^2 \frac{J^3}{X} (\log X)^{o(1)} m(\mathcal{E})^{1/2} \left(\int_{X/2}^{2X}\left| \frac{1}{Q^{3/4}} \sum_{x \leq n < x + Q^{3/4}} \tilde{d}_k(n)\right|^2 dx \right)^{1/2} \\
&\ll \frac{H^2}{\log^{20k \log k} X} \log^{k-1+o(1)} X,
\end{split}
\]
where we have used~\eqref{J-bound}, \eqref{eq:measE''}, and estimated the integral similarly to~(\ref{eq:4thmom}). Thus,
\begin{align*}
A_{i,i'} &\ll \log^{k-1+o(1)} X  \int\int_{|x'-x| \leq H/10} \Biggl( \frac{1}{q_{i,i'}} \sum_{n \in [x,x+H/10]}
|G_{i,x}(n)| |G_{i',x'}(n)| \\
&\qquad + J Q^{-1/4} H \left(\sum_{j=1}^J |g_j(x)|^2\right)^{1/2} \left(\sum_{j=1}^J |g_j(x')|^2\right)^{1/2} \Biggr) dx dx' \\
&\qquad \qquad + \frac{H^2}{\log^{20k \log k} X} \log^{k-1+o(1)} X.
\end{align*}
From Cauchy-Schwarz and the large sieve inequality (see e.g. \cite[formula (7.27)]{ik}) we have
$$
\sum_{n \in [x,x+H/10]} |G_{i,x}(n)| |G_{i',x'}(n)|
\ll H 
\left(\sum_{j \in {\mathcal J}_i} |g_j(x)|^2\right)^{1/2}
\left(\sum_{j \in {\mathcal J}_{i'}} |g_j(x')|^2\right)^{1/2}
$$
and then by a further Cauchy-Schwarz and \eqref{gj2} we conclude that
$$ A_{i,i} \ll H^2 \log^{k-1+o(1)} X  \left( \frac{|{\mathcal J}_i|^{1/2} |{\mathcal J}_{i'}|^{1/2}}{q_{i,i'}} + J^2 Q^{-1/4} + \frac{1}{(\log X)^{20k\log k}} \right).$$
The contribution of the two last terms to~\eqref{eq:Aii'claim} is acceptable from the bounds on $J$ and $Q$, so we are reduced to showig that
$$ \sum_{(i,i') \hbox{ bad}} \frac{|{\mathcal J}_i|^{1/  2} |{\mathcal J}_{i'}|^{1/2}}{q_{i,i'}} \ll J^{5/4}.$$

A similar sum without the sets $\mathcal{J}$ was dealt with in~\cite[(4.4)]{gt-selberg} by Green and Tao and we adapt their argument. The contribution of those $q_{i,i'}$ with $q_{i,i'} \geq J^{3/4}$ is acceptable by the Cauchy-Schwarz inequality, so we may restrict attention to those $(i,i')$ with $q_{i,i'} < J^{3/4}$. By \eqref{ga}, we can bound this contribution by
$$ \sum_{q \leq J^{3/4}} \frac{1}{q} \sum_{a \in \Z/q\Z} \sum_{1 \leq i,i' \leq J'} 
|{\mathcal J}_i|^{1/2} |{\mathcal J}_{i'}|^{1/2} 1_{\| \beta_i - \beta_{i'} - \frac{a}{q} \|_{\T} \leq \frac{3}{Q}}.$$
Let $\Phi: \R \to \R^+$ be an even non-negative Schwartz function with $\Phi(x) \geq 1$ for $x \in [-3,3]$, and whose Fourier transform is supported on $[-1/2,1/2]$.    We may bound the preceding expression by

$$ \sum_{q \leq J^{3/4}} \frac{1}{q} \sum_{a \in \Z/q\Z} \sum_{1 \leq i,i' \leq J'} 
|{\mathcal J}_i|^{1/2} |{\mathcal J}_{i'}|^{1/2} \sum_{m\in \Z} \Phi\left( 
Q\left( \beta_i - \beta_{i'} -\frac{a}{q} + m\right) \right).$$
By the Poisson summation formula, this 
$$ \sum_{q \leq J^{3/4}} \frac{1}{q} \sum_{a \in \Z/q\Z} \sum_{1 \leq i,i' \leq J'} |{\mathcal J}_i|^{1/2} |{\mathcal J}_{i'}|^{1/2} \frac{1}{Q} \sum_{n\in \Z} \hat \Phi\left( \frac{n}{Q} \right) e\left( n\left( \beta_i - \beta_{i'} -\frac{a}{q} \right) \right).$$
Performing the $a$ summation this becomes
$$ \sum_{q \leq J^{3/4}} \sum_{1 \leq i,i' \leq J'} |{\mathcal J}_i|^{1/2} |{\mathcal J}_{i'}|^{1/2} \frac{1}{Q} \sum_{n\in \Z: q|n} \hat \Phi\left( \frac{n}{Q} \right) e\left( n\left( \beta_i - \beta_{i'}  \right) \right).$$
This factorizes as

$$ \frac{1}{Q} \sum_{n \in \Z} \hat \Phi\left(\frac{n}{Q}\right) \left(\sum_{q \leq J^{3/4}: q|n} 1\right) \left|\sum_{i=1}^{J'} |{\mathcal J}_i|^{1/2} e(n\beta_i)\right|^2.$$
From the support of $\hat \Phi$ and H\"older's inequality, we may bound this by
$$ \frac{1}{Q} \left(\sum_{|n| \leq Q} \left(\sum_{q \leq J^{3/4}: q|n} 1\right)^{10}\right)^{1/10}
\left(\sum_{|n| \leq Q} \left|\sum_{i=1}^{J'} |{\mathcal J}_i|^{1/2} e(n\beta_i)\right|^{20/9}\right)^{9/10}.$$
From the large sieve inequality (see e.g. \cite[formula (7.27)]{ik}), the $1/Q$-separated nature of $\beta_i$ and \eqref{J-bound}, \eqref{Q-def} we have
$$\sum_{|n| \leq Q} \left|\sum_{i=1}^{J'} |{\mathcal J}_i|^{1/2} e(n\beta_i)\right|^2 \ll J Q$$
which implies (together with the trivial bound $|\sum_{i=1}^{J'} |{\mathcal J}_i|^{1/2}  e(n\beta_j)| \leq \sum_{i=1}^{J'} |{\mathcal J}_i| = J$) that
$$ \sum_{|n| \leq Q} |\sum_{i=1}^{J'} |{\mathcal J}_i|^{1/2} e(n\beta_i)|^{20/9} \ll J^{11/9} Q$$
and hence the preceding expression may be bounded by
$$ J^{11/10} \left(\frac{1}{Q} \sum_{|n| \leq Q} \left(\sum_{q \leq J^{3/4}: q|n} 1\right)^{10} \right)^{1/10}.$$
But from \eqref{J-bound}, \eqref{Q-def} we have\footnote{One may also invoke truncated divisor sum moment estimates \cite{Bourgain, Ruzsa} here.}
\begin{align*}
\frac{1}{Q} \sum_{|n| \leq Q} \left(\sum_{q \leq J^{3/4}: q|n} 1\right)^{10} 
&\ll \sum_{q_1,\dots,q_{10} \leq J^{3/4}} \frac{1}{Q} \sum_{\substack{|n| \leq Q \\ q_1,\dots,q_{10}|n}} 1 \\
&\ll \sum_{q_1,\dots,q_{10} \leq J^{3/4}} \frac{1}{[q_1,\dots,q_{10}]} \\
&\ll J^{3/4} \sum_{q_1,\dots,q _{10}} \frac{1}{[q_1,\dots,q_{10}]^{11/10}} \\
&\ll J^{3/4} \prod_p \left(1 + O\left( \frac{1}{p^{11/10}}\right )\right) \\
&\ll J^{3/4}
\end{align*}
and the claim follows (with some room to spare).  This concludes the proof of part (i) of Proposition \ref{lve-prop}.

We now briefly discuss the changes needed to handle Proposition \ref{lve-prop}(ii), in which $f \coloneqq \tilde \Lambda 1_{(X,2X]}$.  The main change is to replace the divisor function majorant $\tilde d_k$ by a standard sieve majorant $\nu$ for the (restricted) von Mangoldt function $\tilde \Lambda$, at the level $M$ defined by \eqref{M-def}.  The precise choice of majorant is not of critical importance, but one can for instance use the Goldston-Y{\i}ld{\i}r{\i}m type majorant
\begin{equation}\label{nudef}
 \nu(n) \coloneqq 1_{(X,2X]}(n) \log M \left(\sum_{d|n: d \leq M} \mu(d) \psi\left(\frac{\log d}{\log M}\right)\right)^2
\end{equation}
where $\psi: \R \to \R$ is a smooth function supported on $[-1/2,1/2]$ that equals $1$ at the origin.
From \eqref{M-def} we clearly have the analogue
$$ f(n) \ll \nu(n) \log^{o(1)} X$$
of Lemma \ref{ndic}.  The analogue of Lemma \ref{henriot} for $\nu$ follows from \cite[Theorem D.3]{gt-linear} (see also \cite[Theorem 1.1]{goldston-yildirim-old2}); from this, one can establish Proposition \ref{excep} with $\tilde d_k$ replaced by $\nu$ (and $k$ replaced by $1$).  One can write the majorant $\nu$ in a form similar to \eqref{tdk-def}, but with $M$ replaced by $M^2$ and $d_{k-1}(m)$ replaced by the quantity
$$ \log M \sum_{m_1,m_2: [m_1,m_2] = m} \mu(m_1) \mu(m_2) \log\left(\frac{M}{m_1}\right) \log\left(\frac{M}{m_2}\right).$$
This quantity can be bounded crudely by $O( \log X d_2(m)^2 )$, and one can adapt the proof of Proposition \ref{doc} with $\tilde d_k$ replaced $\nu$ (and $k$ replaced by $1$) without difficulty.  Continuing the remainder of the arguments in this section, we obtain Proposition \ref{lve-prop}(ii).

\end{document}